\documentclass[12pt]{amsart}

\usepackage{amsmath,amsfonts,amssymb}

\theoremstyle{plain}
\newtheorem{lemma}{Lemma}[section]
\newtheorem{theorem}[lemma]{Theorem}

\theoremstyle{definition}

\numberwithin{equation}{section}
\DeclareMathOperator{\var}{var}

\begin{document}

\newcommand{\ZZ}{\mathbb{Z}}
\newcommand{\ZZd}{\mathbb{Z}^{d}}
\newcommand{\RR}{\mathbb{R}}
\newcommand{\RRd}{\mathbb{R}^{d}}
\newcommand{\PP}{\mathbb{P}}
\newcommand{\QQ}{\mathbb{Q}}
\newcommand{\EE}{\mathbb{E}}
\newcommand{\mB}{\mathcal{B}}
\newcommand{\mC}{\mathcal{C}}
\newcommand{\mD}{\mathcal{D}}
\newcommand{\mE}{\mathcal{E}}
\newcommand{\mF}{\mathcal{F}}
\newcommand{\mG}{\mathcal{G}}
\newcommand{\mH}{\mathcal{H}}
\newcommand{\mJ}{\mathcal{J}}
\newcommand{\mL}{\mathcal{L}}
\newcommand{\mkN}{\mathfrak{N}}
\newcommand{\mR}{\mathcal{R}}
\newcommand{\mS}{\mathcal{S}}
\newcommand{\mT}{\mathcal{T}}
\newcommand{\mU}{\mathcal{U}}
\newcommand{\mW}{\mathcal{W}}
\newcommand{\bs}{\backslash}
\newcommand{\half}{\frac{1}{2}}
\newcommand{\fN}{\frac{1}{N}}
\newcommand{\PGibbs}{P^{\beta,u}_{N,\{V_j\}}}
\newcommand{\tZ}{\tilde{Z}}
\newcommand{\tJ}{J_1}
\newcommand{\tK}{\tilde{K}}
\newcommand{\tf}{\tilde{f}}
\newcommand{\tE}{\tilde{E}}
\newcommand{\tL}{\tilde{L}}
\newcommand{\tS}{\tilde{\tau}}
\newcommand{\hd}{\hat{\delta}}
\newcommand{\hP}{\hat{P}}
\newcommand{\hQ}{\hat{Q}}
\newcommand{\hZ}{\hat{Z}}
\newcommand{\hm}{\hat{\mu}}
\newcommand{\ophi}{\overline{\varphi}}
\newcommand{\tphi}{\tilde{\varphi}}
\newcommand{\oF}{\overline{F}}

\title[Disorder and Polymer Pinning]
{The Effect of Disorder on Polymer Depinning Transitions}
\author{Kenneth S. Alexander}
\address{Department of Mathematics KAP 108\\
University of Southern California\\
Los Angeles, CA  90089-2532 USA}
\email{alexandr@usc.edu}
\thanks{Research supported by NSF grant DMS-0405915.}

\date{\today}

\begin{abstract}
We consider a polymer, with monomer locations modeled by the
trajectory of a Markov
chain, in the presence of a potential that
interacts with the polymer when it visits a particular site 0.  We assume that probability of an excursion of length $n$ is given by $n^{-c}\varphi(n)$ for some $1<c<2$ and slowly varying $\varphi$.  Disorder is
introduced by having the interaction vary from one monomer to
another, as a constant $u$ plus i.i.d. mean-0 randomness.  There is a
critical value
of $u$ above which the polymer is pinned, placing a positive fraction
(called the contact fraction) of its monomers
at 0 with high probability.  To see the effect of disorder on the depinning transition,
we compare the contact fraction and free energy (as functions of $u$) to the corresponding annealed system.  We show that for $c>3/2$, at high temperature, the quenched and annealed curves differ significantly only in a very small neighborhood of the critical point---the size of this neighborhood scales as $\beta^{1/(2c-3)}$ where $\beta$ is the inverse temperature.  For $c<3/2$, given $\epsilon>0$, for sufficiently high temperature the quenched and annealed curves are within a factor of $1-\epsilon$ for all $u$ near the critical point; in particular the quenched and annealed critical points are equal.  For $c=3/2$ the regime depends on the slowly varying function $\varphi$.
\end{abstract}
\maketitle

\section{Introduction} \label{S:intro}
We consider the following model for a polymer (or other
one-dimensional object) in a
higher-dimensional space, interacting with a potential located either at a
single site or in a one-dimensional subspace.  There is an
underlying Markov chain $\{ X_i, i \geq 0
\}$, with state space $\Sigma$, representing the ``free'' trajectories of the object in the absence of the
potential.  We assume
the chain is irreducible and aperiodic.
There is a unique site in $\Sigma$ which we
call 0 where the potential can be nonzero; we consider trajectories of length
$N$ starting from state 0 at time 0 and denote the corresponding measure
$P^X_{[0,N]}$.
The potential at 0 at time
$i$ has form $u + V_i$ where the
$V_i$ are i.i.d. with mean zero; we refer to $\{V_i: i \geq 1\}$ as the
\emph{disorder}.  For $N$ and a realization $\{ V_i, i \geq 1
\}$ fixed, we attach a Gibbs weight
\begin{equation} \label{E:GibbsWtQ}
        \exp\left( \beta \sum_{i=1}^N (u + V_i)\delta_{\{x_i = 0 \}} \right)
          P^X_{[0,N]} (x_{[0,N]})
\end{equation}
to each trajectory $x_{[0,N]} = \{ x_i: 0 \leq i \leq N \}$, where $\beta>0$ and $u \in \RR$.  Here $\delta_A$ denotes the
indicator of an event $A$.  This incorporates both the \emph{directed} case in
which we view the object as existing in $\ZZ \times \Sigma$ with configuration
given by the space-time trajectory $(0,0), (1,X_1),..,(N,X_N)$, and the
\emph{undirected} case in which the object exists in $\Sigma$ with
configuration
given by the trajectory $0,X_1,..,X_N$, with $i$ merely an index.

A number of physically disparate phenomena are subsumed in this formulation,
which was studied in \cite{AS06}, \cite{Gi06} and \cite{GT06} in the mathematics literature; closely related models have been studied extensively in the physics literature, for example in \cite{DHV92}, \cite{GG96}, \cite{MB93}, \cite{NN01} and \cite{NZ95}.  Taking $\Sigma = \ZZd$ we obtain a directed
quenched random copolymer in
$d+1$ dimensions interacting, with strength that depends on the random monomer
types, with a potential located in the first coordinate axis.
Here the $(i+1)$st monomer is located at
$(i,X_i)$ and the realization $\{V_i\}$, representing the random sequence of monomer types, is fixed.  Alternatively, we may have a
uniform polymer but spatial randomness in the potential (as when the
``polymer'' is actually a flux tube in a superconducting material with
one-dimensional defects--see
\cite{GT01}, \cite{NV93}.)  When $d=1$ and $X_i \geq 0$, the state $X_i$ may represent
the height at location $i$ of an interface above a wall in two
dimensions, which attracts or repels the interface with randomly spatially varying
strength--see \cite{DHV92}, \cite{FLNO88}.

For consistency, here we always view the index $i$ as time.

The main questions of physical interest are (i) whether for given $\beta,u$ 
the polymer is ``pinned'', meaning roughly that it places a positive fraction
of its monomers at 0 for large $n$ with high probability; (ii) the location and nature of the depinning transition, if any, as we vary $\beta$ and/or $u$; and (iii) the effect of the disorder,
as seen by comparing the transition to the annealed case (which is effectively the same as $V_i 
\equiv$
constant.)

The partition function and
finite-volume
Gibbs distribution on length-$N$ trajectories, corresponding to
\eqref{E:GibbsWtQ}, are denoted
$Z^{\{V_i\}}_{[0,N]}(\beta,u)$ and
$\mu^{\beta,u,\{V_i\}}_{[0,N]}$, respectively; the corresponding expectation is denoted $\langle \cdot \rangle^{\beta,u,\{V_i\}}_{[0,N]}$.  We omit the $\{V_i\}$ when $V_i
\equiv 0$.  We write $P^V_{[a,b]}$ for the distribution of $(V_a,..,V_b)$ and $P^V$ 
for  $P^V_{[0,\infty)}$. (Here
and wherever we deal with indices, we take intervals to mean intervals of
integers.)  The corresponding expectation and conditional expectation are
denoted $\langle \cdot \rangle^V_{[a,b]}$ and $E^V_{[a,b]}(\cdot \mid 
\cdot)$.  Let
\[
        L_N = \sum_{i=1}^N \delta_{\{X_i = 0 \}}.
\]
For fixed $\beta, u$ we say the polymer is \emph{pinned} at
$(\beta,u)$ if for some
$\delta > 0$,
\[
        \lim_N \mu^{\beta,u,\{V_i\}}_{[0,N]}( L_N > \delta N) = 1, \quad
          P^V-\text{a.s.}
\]
There is a (possibly infinite)
critical $u_c(\beta,\{V_i\})$ such that the polymer is pinned
for $u > u_c(\beta,\{V_i\})$ and not pinned for $u < u_c(\beta,\{V_i\})$.
In \cite{AS06} it was established that self-averaging  holds in the
sense that there
is a nonrandom  \emph{quenched critical point}
$u_c^q = u_c^q(\beta)$ such that $u_c(\beta,\{V_i\}) = u_c^q(\beta)$ with
$P^V$-probability one.  The
\emph{deterministic critical
point} $u_c^d = u_c^d(\beta)$ is the critical point for the \emph{deterministic
model}, which is the case
$V_i \equiv 0$.  The \emph{annealed model} is obtained (provided the moment
generating function $M_V(\beta) =
\langle e^{\beta V_1} \rangle^V$ is finite) by averaging the Gibbs
weight \eqref{E:GibbsWtQ} over the disorder; the annealed model at
$(\beta,u)$ is
thus the same as the deterministic model at $(\beta,u + \beta^{-1}\log
M_V(\beta))$, and the corresponding
\emph{annealed critical point} is $u_c^a = u_c^a(\beta) = u_c^d(\beta) -
\beta^{-1}\log M_V(\beta)$.  It is not hard to show that
\begin{equation} \label{E:critpts}
        u_c^a \leq u_c^q \leq u_c^d;
\end{equation}
see \cite{AS06}.  For the annealed case, we denote the partition function and finite-volume Gibbs distribution by $Z_{[0,N]}(\beta,u)$ and $\mu_{[0,N]}^{\beta,u}$, respectively, and for the deterministic case we denote them by $Z^0_{[0,N]}(\beta,u)$ and $\mu_{[0,N]}^{\beta,u,0}$, respectively.
The free energies for the deterministic, annealed
and quenched
models are given by
\[
    f^d(\beta,u) = \frac{1}{\beta} \lim_N \fN \log Z^0_{[0,N]}(\beta,u),
\]
\[
    f^a(\beta,u) = \frac{1}{\beta} \lim_N \fN \log Z_{[0,N]}(\beta,u) = f^d(\beta,u+\beta^{-1}\log M_V(\beta)),
\]
\[
    f^q(\beta,u) = \frac{1}{\beta} \lim_N \fN \log Z^{\{V_i\}}_{[0,N]}(\beta,u).
\]
The $P^V$-a.s. existence and the nonrandomness of the last limit are proved in \cite{GT06}, for the cases we consider here.
Clearly $f^d(\beta,u)$ depends only on $\beta u$, and
$f^a(\beta,u)$ depends only on $\beta u + \log M_V(\beta)$.  The \emph{specific
heat exponent} $\alpha_d$ in the deterministic case is given by
\[
    2 - \alpha_d = \lim_{\Delta \searrow 0}
      \frac{\log f^d(\beta,u_c^d(\beta) + \Delta)}{\log \Delta};
\]
clearly we get the same annealed exponent $\alpha_a$ if we replace
$d$ with $a$ in
this definition.  For the quenched case the same definition applies with $d$
replaced by $q$, provided the limit exists.

A related quantity of interest is the \emph{contact
fraction}, defined in the deterministic model to be the value $C =
C^d(\beta,u)$
for which
\begin{equation} \label{E:contfrac}
     \lim_{\epsilon \searrow 0} \lim_N \nu^{\beta,u}_{[0,N]}\left( \frac{L_N}{N} \in
      (C-\epsilon,C+\epsilon) \right) = 1;
\end{equation}
the existence of such a $C$ is established in \cite{AS06}.  In the
annealed case
the contact fraction is
\[
    C^a(\beta,u) = C^d(\beta,u+\beta^{-1}\log M_V(\beta)).
\]
Now $f^d$ and $f^a$ are convex
functions of $u$, and we have by standard methods that
\[
    C^d(\beta,u) = \frac{\partial}{\partial u} f^d(\beta,u), \quad
      C^a(\beta,u) = \frac{\partial}{\partial u} f^a(\beta,u)
\]
for all non-critical $u$; the necessary differentiability of $f^a(\beta,\cdot)$ follows readily from the differentiability and strict convexity established in the proof of Lemma 2.2 of \cite{AS06} for the function $F$ defined in \eqref{E:Fdef} below.  In the quenched case, differentiability of $f^q(\beta,\cdot)$ is proved in \cite{GT06b} when the underlying Markov chain is simple random walk on $\ZZ$, but is not known in general, forcing us to define $C^{q,-}(\beta,u)$ and $C^{q,+}(\beta,u)$ to be the left and right derivatives respectively of the convex function $f^q(\beta,\cdot)$ at $u$, and then define 
\[
  \mD(\beta) = \{u \in \RR: u_c^q(\beta): C^{q,+}(\beta,u) = C^{q,-}(\beta,u)\},
  \]
a set for which the complement is at most countable.  For $u \in \mD(\beta)$ we denote the common value $\frac{\partial}{\partial u} f^q(\beta,u)$ by $C^q(\beta,u)$.  
From convexity we have for fixed $\beta$ that 
\[
  \left\langle \frac{L_N}{N} \right\rangle^{\beta,u,\{V_i\}}_{[0,N]} = \frac{1}{\beta} \frac{\partial}{\partial u} \left(
    \fN \log Z^{\{V_i\}}_{[0,N]}(\beta,u) \right) \to \frac{\partial}{\partial u} 
    f^q(\beta,u) \text{ for all } u \in \mD(\beta)
  \]
and 
\[
  C^{q,-}(\beta,u) \leq \liminf_N \left\langle \frac{L_N}{N} \right\rangle^{\beta,u,\{V_i\}}_{[0,N]}
    \leq \limsup_N  \left\langle \frac{L_N}{N} \right\rangle^{\beta,u,\{V_i\}}_{[0,N]} \leq C^{q,+}(\beta,u) 
    \text{ for all } u \notin \mD(\beta), 
    \]
both $P^V$-a.s.

The case most commonly considered in the literature has $\Sigma = \ZZd$ and
$\{X_i\}$ a symmetric simple random walk.  In keeping with our requirement of
aperiodicity we modify this by considering the case in which $X_i$ is
the location
of the walk at time $2i$; we call this the \emph{symmetric simple
random walk case in d+1 dimensions}.  Let $\tau_i$ denote the time of the $i$th return to 0 by the chain,
with $\tau_0 = 0$, and let $E_i = \tau_i - \tau_{i-1}$ denote the $i$th excursion
length.  Our interest here is in general Markov chains
$\{X_i\}$ which satisfy
\begin{equation} \label{E:tails}
      P^X(E_1 = n) = n^{-c}\varphi(n)
\end{equation}
for some $c \in (1,2)$ and slowly varying function $\varphi$ on 
$[1,\infty)$.  This includes the symmetric simple random walk case in one
and (by virtue of Theorem \ref{T:transient}) three
dimensions, where $c
=  3/2$ and $\varphi(n)
\sim K$ for some $K>0$.  We do not, however, include the cases of $c=1$ and $c\geq 2$, because the technical details, and many of the heuristics, are quite different in these cases, as we will discuss further in Section \ref{S:excludedc}.  We focus instead on our main purpose which is to understand the role of the tail exponent $c$ in the effect of quenched disorder; as we will see, the main distinction is between $c<3/2$ and $c>3/2$. Among excluded examples are the symmetric simple random walk case in $2+1$ dimensions (where $c=1$ and $\varphi(n) \sim K(\log n)^{-2}$), and in $d+1$ dimensions for $d \geq 4$ (where $c=d/2 \geq 2$.)
For recurrent chains satisfying
\eqref{E:tails} it is easily seen (see \cite{AS06}) that $u_c^d(\beta) = 0$ for
all
$\beta>0$, and hence
\begin{equation} \label{E:uca}
     u_c^a(\beta) = -\beta^{-1}\log M_V(\beta),
\end{equation}
and, again from \cite{AS06}, for the
deterministic or annealed model the transition is first order if and only if
$\langle E_1 \rangle^X < \infty$; in particular it is first order for $c>2$ 
but not for $c<2$.  It is proved in \cite{GT06} 
that (again as known nonrigorously from the physics literature--see e.g. \cite{Fi84})
the annealed specific heat exponent is $(2c-3)/(c-1)$.  Based on the physics
literature (\cite{DHV92}, \cite{Fi84}), the following is believed for the non-first-order cases:
\begin{itemize}
\item[(i)] for $3/2 < c < 2$ (positive annealed specific heat exponent) the
depinning transition is altered by the disorder;
\item[(ii)] for $c < 3/2$ (negative annealed specific heat exponent)
the depinning
transition is not altered by the disorder.
\end{itemize}
Here physicists generally take ``altered'' to mean that the specific heat
exponent is different, but a disorder-induced change in the critical
point is also of interest (\cite{CGG06}, \cite{DHV92}, \cite{Mo00}, \cite{NN01}).  This produces the
question, ``change relative to what, $u_c^a$ or $u_x^q$?''  It appears most
natural to ask whether the factor $e^{\beta(u+V_i)}$ in
\eqref{E:GibbsWtQ} gives the same critical $u$ as if it were replaced by its mean
$e^{\beta u} \langle e^{\beta V_i} \rangle^V$, which is equivalent to asking
whether $u_c^q(\beta) = u_c^a(\beta)$.  The question is of interest in part because it is intertwined with questions of just how the polymer depins as the quenched critical point is approached, or, put differently, of what ``strategy'' the polymer uses to stay pinned when near the critical point--see \cite{CGG06}.  For example, do long stretches depin, leaving the polymer attached only where the disorder is exceptionally favorable, or is the depinning more uniform?  Some of our results may be interpreted as saying that for $1<c < 3/2$ the depinning is quite uniform in the quenched system all the way to the critical point, and for $3/2<c<2$ the depinning is quite uniform at least until very close to the critical point, at least for high temperatures.

It was proved in \cite{GT06} that disorder does alter the
critical behavior when $3/2 < c < 2$, in that the quenched specific 
heat exponent (assuming it exists)
becomes non-positive, i.e. the free energy increases no more than quadratically
as $u$ increases from $u_c^q$.  This brings us to one of the main questions we consider here:  just how are the free energy and contact fraction curves (as functions of $u$) altered by the disorder?  When the disorder changes the specific heat exponent from an annealed value $\alpha_a > 0$ to a quenched value $\alpha_q \leq 0$, does a large section of the contact fraction curve $\Delta \mapsto C^a(\beta,u_c^a(\beta) + \Delta) \sim \Delta^{1-\alpha_a}$ change to approximate $\Delta \mapsto \Delta^{1-\alpha_q}$ instead, or does significant change only occur very close to $\Delta=0$?  We will show that the latter is the case, at least for small $\beta$.

Regarding the possible difference in critical points between quenched and annealed systems, since
$u_c^d(\beta) = 0$,
it is reasonable to ask how close
$u_c^q(\beta)/u_c^a(\beta)$ is to 1.  As discussed in \cite{AS06},
this is related
to the following ``sampling'' point of view.  We may think of the
Markov chain as
choosing a sample from the realization $\{V_1,..,V_n\}$ through the
times of its
returns to 0.  This sample is of course not i.i.d., and we expect that the
probabilities for large deviations of the average of the sampled
$V_i$'s will be
smaller than the corresponding probabilities for an i.i.d. sample.
Roughly, the
more this sampling procedure is ``efficient'' in the sense that certain
large-deviation probabilities are not too different from the i.i.d. case, the
closer $u_c^q(\beta)/u_c^a(\beta)$ is to 1.  For small $\beta$, the size of the
relevant large deviations is of order
$\beta$ (see \cite{AS06}.)

The analogous question has been considered for a related model, a
polymer at a
selective interface in $1+1$ dimensions.  Here the horizontal axis
represents an
interface separating two solvents which differentially attract or repel each
monomer, with $u+V_i$ representing roughly the preference of monomer
$i$ for the
solvent above the interface.  The Markov chain $\{X_i\}$ is symmetric simple
random walk on $\ZZ$, so the excursion tail exponent is $c=3/2$.  The
main mathematical
difference from the model here is that the factor
$\delta_{\{x_i = 0 \}}$ in \eqref{E:GibbsWtQ} is replaced by
$-\delta_{\{x_i < 0 \}}$.  It has been conjectured in the physics \cite{Mo00} and mathematics
(\cite{BG04}, \cite{Gi06}) literatures that for this model,
\begin{equation} \label{E:selective}
     \lim_{\beta \to 0} \frac{u_c^q(\beta)}{u_c^a(\beta)} < 1.
\end{equation}
The limit being less than 1 says that the above-discussed sampling
procedure does
not become fully efficient even in the high-temperature limit, where
the size of
the relevant large deviation approaches 0.

For technical convenience and clarity of exposition, we will restrict attention here to Gaussian
disorder, but our results are easily extended to other distributions with a finite
exponential moment.  Our results imply that, in contrast to
\eqref{E:selective}, for the
model \eqref{E:GibbsWtQ} with
$1<c<2$, the limit in \eqref{E:selective} is 1; in fact for $c<3/2$ the ratio is equal to 1 at least for sufficiently small $\beta$.  We do not view this as
casting any doubt on the conjecture \eqref{E:selective} for the selective
interface model, however, for the following reason:  the factor
$\delta_{\{x_i < 0
\}}$ in the Gibbs weight means that the sample selected by the Markov chain
from $\{V_1,..,V_n\}$ consists of blocks $V_j,V_{j+1},..,V_{j+k}$ of
consecutive
disorder values, corresponding to excursions of $\{X_i\}$ below the axis.  This
``block'' aspect may make the sampling procedure inherently less efficient (for
large deviations) in the selective interface model, compared to our model
\eqref{E:GibbsWtQ}.

  From \eqref{E:uca}, in the Gaussian case we have $u_c^a(\beta) = -\beta/2$ for
all $\beta>0$.

Associated to a slowly varying function $\eta$ there is a conjugate slowly varying function 
$\eta^*$ characterized (up to asymptotic equivalence) by the property that 
\[
  \eta^*(x\eta(x)) \sim \frac{1}{\eta(x)} \quad \text{as } x \to \infty;
  \]
see \cite{Se76}.  Here ``$\sim$'' means the ratio converges to 1.
For many common slowly varying functions such as 
$\eta(x) = (\log x)^a$ with $a \in \RR$, 
we have $\eta^* = 1/\eta$, but this is not the case for some functions which are
``barely slowly varying'' such as $\eta(x) = x^{1/\log \log x}$.  Given $\varphi$ as in
\eqref{E:tails} and $a>0$, define
\[
  \ophi_a(x) = \frac{1}{\varphi(x^{1/a})}, \qquad \hat{\varphi}_a(x) = \ophi_a^*(x)^{-1/a},
    \qquad \tphi(x) = \sum_{n \leq x} n^{-1}\varphi(n)^{-2},
  \]
which are all slowly varying.

Our first result is for the annealed case, included mainly so that the quenched
case can be compared to it.  Most of this result appeared in \cite{GT06}; we have
made only minor additions.  Throughout the paper,
$K_1,K_2,...$ are constants which depend only on the distributions of $V_1$ and
$\{X_i\}$, unless otherwise specified. 

\begin{theorem} \label{T:anneal}
Suppose that $\{V_i, i \geq 1\}$ are i.i.d. standard Gaussian random
variables, and the Markov chain $\{X_i\}$ is recurrent, satisfying
\eqref {E:tails} with $1<c<2$.  Then $u_c^a(\beta) = -\beta/2$ for all
$\beta>0$,
and there exist constants $K_i$, depending only on
$c, \varphi$, such that
\[
    \beta f^a\left(\beta,-\frac{\beta}{2} + \Delta\right) \sim
      K_1 (\beta\Delta)^{1/(c-1)} \hat{\varphi}_{c-1} \left( \frac{1}{\beta\Delta}
      \right) \quad \text{as } \beta\Delta \to 0,
\]
and
\[
    C^a\left(\beta,-\frac{\beta}{2} + \Delta\right) \sim
      K_2 (\beta\Delta)^{(2-c)/(c-1)} \hat{\varphi}_{c-1} \left( \frac{1}{\beta\Delta}
      \right) \quad \text{as } \beta\Delta \to 0.
\]
In particular, the annealed specific heat exponent is $(2c-3)/(c-1)$.
\end{theorem}

Note that writing the parameter $u$ as $-\frac{\beta}{2} + \Delta$
in Theorem
\ref{T:anneal} makes the annealed free energy and contact fraction functions of
$\beta\Delta$ only, as is reflected in the asymptotic approximations given
there.

The next theorem confirms the prediction that disorder does not alter the
critical behavior when $c < 3/2$, at high temperatures; we can even include
certain cases when
$c = 3/2$, though not the symmetric simple random walk case in 1+1 dimensions.

\begin{theorem} \label{T:smallc}
Suppose that $\{V_i, i \geq 1\}$ are i.i.d. standard Gaussian random
variables, and
the Markov chain $\{X_i\}$ is recurrent, satisfying \eqref{E:tails}
with either $1 < c < 3/2$, or $c= 3/2$ and $\sum_n n^{-1} \varphi(n)^{-2} <
\infty$.  Then given $\epsilon>0$, provided $\beta>0$ and $\beta\Delta>0$ are 
sufficiently small we have 
\begin{equation} \label{E:highfree1}
    1 - \epsilon \leq \frac{
    f^q\left(\beta,-\frac{\beta}{2} + \Delta\right)}{
      f^a\left(\beta,-\frac{\beta}{2} + \Delta\right)} \leq 1,
      \quad \left| \frac{
      C^{q,\pm}\left(\beta,-\frac{\beta}{2} + \Delta\right)}{
      C^a\left(\beta,-\frac{\beta}{2} + \Delta\right)} - 1 \right| \leq \epsilon,
\end{equation}
so that in particular $u_c^q(\beta) = u_c^a(\beta) = -\beta/2$.
\end{theorem}

As we will see, when $c=3/2$, $\tphi(n)$ is proportional to the mean overlap under $P^X$ between 
$\{X_i\}$ and an independent copy $\{Y_i\}$ of the chain over length $n$, that is, the mean 
number of $i \leq n$ with $(X_i,Y_i) = (0,0)$.  The condition in Theorem \ref{T:smallc} that $\sum_n
n^{-1} \varphi(n)^{-2} < \infty$, i.e. that $\tphi$ is bounded, is thus equivalent (for $c=3/2$) to the condition that
$(X_i,Y_i)$ is transient.

The next theorem quantifies the change in the critical curve caused by the presence of disorder, for $3/2<c<2$.  It
shows that at high
temperatures, the significant alteration is confined to a very small interval
above
$u_c^a$, with length of order $\beta^{1/(2c-3)} \hat{\varphi}_{c-\frac{3}{2}}(\beta^{-1})^{1/2}$.  In particular, the critical points
$u_c^q$ and $u_c^a$ differ by no more than $\beta^{1/(2c-3)} 
\hat{\varphi}_{c-\frac{3}{2}}(\beta^{-1})^{1/2}$.  Note the exponent $1/(2c-3)$ here is greater than 1.

\begin{theorem} \label{T:bigc}
Suppose that $\{V_i, i \geq 1\}$ are i.i.d. standard Gaussian random
variables, and
the Markov chain $\{X_i\}$ is recurrent, satisfying \eqref{E:tails}
for some $3/2 < c < 2$.  Then 
there exists $K_{3}$ as follows.  Let
\[
  \Delta_0 = \Delta_0(\beta) = 
    K_{3} \beta^{1/(2c-3)} \hat{\varphi}_{c-\frac{3}{2}}\left( \frac{1}{\beta} \right)^{1/2}.
  \]
Given $\epsilon>0$, there exists $K_{4}(\epsilon)$ such that 
for all sufficiently small $\beta>0$ and $\beta\Delta>0$ we have
\begin{equation} \label{E:lowfree1}
    \beta f^q\left(\beta,-\frac{\beta}{2}+\Delta\right) \leq \frac{\Delta^2}{2},
      \quad C^{q,\pm}\left(\beta,-\frac{\beta}{2} + \Delta\right) \leq
      \frac{2}{\beta} \Delta
      \quad \text{if } \Delta \leq \Delta_0,
\end{equation}
\begin{equation} \label{E:highfree2}
    1 - \epsilon \leq \frac{
    f^q\left(\beta,-\frac{\beta}{2} + \Delta\right)}{
      f^a\left(\beta,-\frac{\beta}{2} + \Delta\right)} \leq 1,
      \quad \left| \frac{
      C^{q,\pm}\left(\beta,-\frac{\beta}{2} + \Delta\right)}{
      C^a\left(\beta,-\frac{\beta}{2} + \Delta\right)} - 1 \right| \leq \epsilon
      \quad \text{if } \Delta \geq K_{4}\Delta_0.
\end{equation}
Consequently we have
\begin{equation}
    -\frac{\beta}{2} \leq u_c^q(\beta) \leq
    -\frac{\beta}{2} + K_{4}\Delta_0(\beta),
\end{equation}
and therefore
\begin{equation} \label{E:ratio}
     \lim_{\beta \to 0} \frac{u_c^q(\beta)}{u_c^a(\beta)} = 1.
\end{equation}
The critical behavior differs from the annealed case in that either
$u_c^q(\beta)
\neq -\beta/2$ or the specific heat exponent, if it exists, is not
$(3-2c)/(c-1)$.
\end{theorem}

The free-energy and contact-fraction bounds in \eqref{E:lowfree1} are actually
valid for all $\Delta>0$, but when $\Delta > \Delta_0(\beta)$
the free-energy upper bound $\Delta^2/2$ exceeds the trivial upper bound 
$\beta f^a\left(\beta,-\frac{\beta}{2} + \Delta\right)$, as we will obtain from Theorem \ref{T:anneal}, so 
it provides no useful information.  $\Delta_0(\beta)$ is also (up to a constant) the magnitude of $\Delta$ for which the upper bound $2\Delta/\beta$ is equal to the annealed contact fraction, as we will show.  We will also see later that a value $M = M(\beta\Delta)$ of order $1/\beta f^a(\beta,u)$ makes a natural (annealed) correlation length for the problem.  For a block of length $M$ we may view the average value of $u+V_i$ in the block as a sort of ``effective $u$'' for that block; this effective $u$ fluctuates by order $M^{-1/2}$ from block to block.  $\Delta_0$, as we will show, can be characterized (up to a constant) as the value such that for $\Delta \ll \Delta_0$ the fluctuations of the effective $u$ are $M^{-1/2} \gg \Delta$, while for $\Delta \gg \Delta_0$ we have $M^{-1/2} \ll \Delta$.  Thus for $\Delta \leq \Delta_0$ a substantial fraction of blocks have an effective $u$ below $u_c^a$, but this does not occur for $\Delta \gg \Delta_0$.

In the borderline case $c=3/2$, where the annealed specific heat exponent is 0,
outside of the case covered by Theorem \ref{T:smallc} we are unable to say
whether the critical behavior is altered by the disorder; there is disagreement on this question even in the physics literature (\cite{DHV92}, \cite{FLNO88}).  However, in this
case the free energy changes significantly in at most an even smaller
interval above $u_c^a$, with length of order $o(\beta^r)$ for \emph{all} $r>1$, in fact
$O(e^{-K/\beta^2})$ for some constant $K$ in the case of simple random walk in $1+1$ or $3+1$ dimensions.
This may help explain why nonrigorous techniques such as renormalization group
methods and numerical simulation do not provide consistent predictions for the
$c=3/2$ case.

\begin{theorem} \label{T:c32}
Suppose that $\{V_i, i \geq 1\}$ are i.i.d. standard Gaussian random
variables, and
the Markov chain $\{X_i\}$ is recurrent, satisfying \eqref{E:tails}
with $c=3/2$ and $\sum_n n^{-1} \varphi(n)^{-2} = \infty$.
Then given $0<\epsilon<1$, there exists $K_5 = K_5(\epsilon)$
as follows.  Let
\begin{equation} \label{E:Delta0def}
  \Delta_0 = \Delta_0(\beta,\epsilon) = 
    \frac{ K_{5} \varphi\left(\tphi^{-1}\left( \frac{K_{6}}{\beta^2} \right) \right) }
    { \tphi^{-1}\left( \frac{K_{6}}{\beta^2} \right)^{1/2} }.
    \end{equation}
Provided $\beta$ and $\beta\Delta$ are sufficiently small and $\Delta \geq \Delta_0$, we have
\[
    1 - \epsilon \leq \frac{
    f^q\left(\beta,-\frac{\beta}{2} + \Delta\right)}{
      f^a\left(\beta,-\frac{\beta}{2} + \Delta\right)} \leq 1,
      \quad \left| \frac{
      C^{q,\pm}\left(\beta,-\frac{\beta}{2} + \Delta\right)}{
      C^a\left(\beta,-\frac{\beta}{2} + \Delta\right)} - 1 \right| \leq \epsilon.
\]
In particular we have
\[
    -\frac{\beta}{2} \leq u_c^q(\beta) \leq
    -\frac{\beta}{2} + \Delta_0(\beta,\epsilon),
\]
and therefore
\begin{equation} \label{E:ratio2}
     \lim_{\beta \to 0} \frac{u_c^q(\beta)}{u_c^a(\beta)} = 1.
\end{equation}
\end{theorem}

As we will see, the condition $\Delta \geq \Delta_0$ says that the mean overlap mentioned above is at most of order $1/\beta^2$, over one correlation length $M$.

In the symmetric simple random walk case in $1+1$ or $3+1$ dimensions, we have $c=3/2$ and  
$\varphi$ asymptotically constant, so 
\[
  \tphi(x) \sim \log x \quad \text{and} \quad \Delta_0 \leq K_{7} e^{-K_{8}/\beta^2}.
  \]
  
For a transient chain satisfying \eqref{E:tails}, it is proved in \cite{AS06} that 
\begin{equation} \label{E:ucd}
  u_c^d(\beta) = - \beta^{-1} \log P^X(E_1 < \infty).
  \end{equation}
To make Theorems \ref{T:smallc}--\ref{T:c32} useful in the transient case, we need a stronger result.  Given a measure $P^X$ which makes $\{X_i\}$ a transient Markov chain, we define a modified measure $P^{X,R}$ by 
\[
  P^{X,R}(E_1 = n) = P^X(E_1 = n \mid E_1 < \infty),
  \]
keeping the independence of the $E_i$'s and the same conditional distribution given the $E_i$'s.  It is easily checked that the resulting ``recurrent-ized'' process is still a Markov chain, so Theorems \ref{T:smallc}--\ref{T:c32} apply to it.  We denote the corresponding free energy and contact fraction by $f^q_R(\beta,u)$ and $C^{q,\pm}_R(\beta,u)$ (or $C^{q,\pm}_R(\beta,u)$ at points of differentiability), respectively.

\begin{theorem} \label{T:transient}
Suppose that $\{V_i, i \geq 1\}$ are i.i.d. standard Gaussian random
variables, and
the Markov chain $\{X_i\}$ is transient, satisfying \eqref{E:tails}
for some $1 \leq c < 2$.  Then
\[
  f^q(\beta,u) = f^q_R(\beta,u - u_c^d(\beta)), \qquad C^{q,\pm}(\beta,u) = C^{q,\pm}_R(\beta,u - u_c^d(\beta)) 
    \quad \text{for all } \beta>0, u \in \RR.
  \]
\end{theorem} 

\section{Preliminaries on Large Deviations and Asymptotics} \label{S:prelims}
In this section we summarize various basic results on large deviations specialized to our context, and certain asymptotics including Theorem \ref{T:anneal}, for ready reference later.  
We assume throughout that \eqref{E:tails} holds with $1<c<2$.
Let $M_E$ denote the moment generating function of $E_1$ and
\begin{equation} \label{E:IEdef}
    I_E(t) = \sup_{x \leq 0}(tx - \log M_E(x))
\end{equation}
the large-deviation rate function of $E_1$.
Fix $\beta > 0$ and $u > -\frac{\beta}{2}$, and define $\Delta$ by $u =
-\frac{\beta}{2} + \Delta$.  From \cite{AS06}, we have
\[
    \lim_{N \to \infty} \fN \log P^X(L_N \geq \delta N) = -\delta
I_E(\delta^{-1}),
\]
and the free energy of the deterministic model is given by
\[
    \beta f^d(\beta,u) = \lim_N \frac{1}{N} \log Z_{[0,N]}(\beta,u)
    = \frac{1}{\beta} \sup_{\delta \in (0,1)} (\beta u\delta - \delta
    I_E(\delta^{-1})).
\]
Hence the annealed free energy is given by
\begin{align} \label{E:annfreedef}
    \beta f^a(\beta,u) &= \beta f^d(\beta,u+\beta^{-1}\log M_V(\beta)) \\
    &= \sup_{\delta \in (0,1)} \left( \delta \left( \beta u
      + \log M_V(\beta) \right) - \delta I_E(\delta^{-1}) \right) \notag \\
    &= \sup_{\delta \in (0,1)} \left( \beta\Delta\delta - \delta
      I_E(\delta^{-1}) \right). \notag
\end{align}
Let $\delta^* = \delta^*(\beta\Delta)$ be the unique value of $\delta$ which
achieves the supremum in
\eqref{E:annfreedef}; it is easily shown that $\delta^* = C^a(\beta,u)$.
The uniqueness here follows from strict convexity of $\delta I_E(\delta^{-1})$;
see \cite{AS06}.
Let $\delta_n = E^X(L_n/n)$.  From Lemma \ref{L:Erickson} below (applied with $\varphi(n)/(c-1)$ in place of $\varphi(n)$) we see that
\begin{equation} \label{E:deltan}
   \delta_n \sim \frac{(c-1)\Gamma(2-c)}{\Gamma(c-1)}n^{-(2-c)}\varphi(n)^{-1}.
\end{equation}
Let
\begin{equation} \label{E:Mdef}
   M = M(\beta\Delta) = \min\{n: \delta_n \leq \delta^*\},
\end{equation}
and observe that $\delta^* \geq \delta_M \geq (1 - M^{-1})\delta^*$.   $M$ 
serves as a
correlation length for the pinned polymer--we will see that the free energy gain from pinning is of order 1 over length $M$.  It can also be shown that excursions much
longer than $M$ are rare, and on length scales shorter than $M$, the pinned Markov chain in many senses ``looks like'' the underlying ``free'' chain with law $P^X$, but we will not formalize or prove these statements here.

Since $\langle E_1 \rangle^X = \infty$, there exists $\alpha_0 =
\alpha_0(\beta\Delta)>0$ such that
\begin{equation} \label{E:alphadef}
   (\log M_E)'(-\alpha_0) =
   \frac{\langle E_1e^{-\alpha_0 E_1} \rangle^X}{\langle e^{-\alpha_0 E_1} \rangle^X} =
    \frac{1}{\delta^*}.
\end{equation}
From basic large-deviations theory, as shown in 
\cite{GT06} we have
\begin{equation} \label{E:LDfacts}
   \alpha_0 =
     \beta f^a(\beta,u) \quad \text{and } \quad \beta\Delta = -\log M_E(-\alpha_0).
\end{equation}
By \eqref{E:deltan},
\begin{equation} \label{E:deriv1}
   (\log M_E)'(-\alpha_0) = \frac{1}{\delta^*} \sim \frac{1}{\delta_M} \sim
     \frac{\Gamma(c-1)}{(c-1)\Gamma(2-c)} M^{2-c}\varphi(M)
     \quad \text{ as } \beta\Delta \to 0.
\end{equation}
A routine calculation shows that the derivatives of $\log M_E$ satisfy
\begin{equation} \label{E:logMEderivs}
   (\log M_E)^{(k)}(-x) \sim M_E^{(k)}(-x) \sim \Gamma(k+1-c)x^{-(k+1-c)} 
\varphi(x^{-1}) \quad
     \text{as } -x \nearrow 0, \quad k \geq 1,
\end{equation}
so for every fixed $\nu>0$,
\begin{equation} \label{E:deriv2}
   (\log M_E)'\left( -\frac{\nu}{M} \right) \sim 
\frac{\Gamma(2-c)M^{2-c}\varphi(M)}
     {\nu^{2-c}} \quad \text{as } \beta\Delta \to 0 \text{ (equivalently, as }
     M \to \infty.)
\end{equation}
Taking $\nu$ to be $K_{9}$ given by
\[
   \frac{\Gamma(2-c)}{K_{9}^{2-c}} = \frac{\Gamma(c-1)}{(c-1)\Gamma(2-c)},
\]
we see from \eqref{E:deriv1} and \eqref{E:deriv2} that
\begin{equation} \label{E:ma0}
   M\beta f^a(\beta,u) = M\alpha_0 \to K_{9} \quad \text{as } 
\beta\Delta \to 0.
\end{equation}

For $\delta>0$ let
\[
   x_0 = x_0(\delta) = ((\log M_E)')^{-1}(\delta^{-1}),
   \]
that is, $x_0<0$ is the point where the sup in \eqref{E:IEdef} is achieved for
$t=\delta^{-1}$.  Note that $x_0(\delta^*) = -\alpha_0$.
 From \eqref{E:logMEderivs} we have as $\delta \to 0$
\begin{equation} \label{E:deltavsx0}
   \frac{1}{\delta} = (\log M_E)'(x_0) \sim \frac{\Gamma(2-c)}{|x_0|^{2-c}} \varphi\left( \frac{1}{|x_0|} \right)
     = \frac{\Gamma(2-c)}{|x_0|^{2-c}} 
     \overline{ \left( \frac{1}{\varphi} \right) }_{2-c}\left(  \frac{1}{|x_0|^{2-c}} \right),
   \end{equation}
so
\[
  \frac{1}{\delta} \overline{ \left( \frac{1}{\varphi} \right) }^*_{2-c}\left(  \frac{1}{\delta} \right)
    \sim \frac{\Gamma(2-c)}{|x_0|^{2-c}}
  \]
so
\begin{equation} \label{E:x0vsdelta}
   x_0(\delta) \sim -K_{10} \delta^{1/(2-c)} 
     \widehat{ \left( \frac{1}{\varphi} \right) }_{2-c}\left( \frac{1}{\delta} \right).
\end{equation}
Letting
\begin{equation} \label{E:Fdef}
   F(\delta) = \beta\Delta \delta - \delta I_E(\delta^{-1})
   \end{equation}
we have 
\[
   -\delta^2 F''(\delta) = \delta^{-1}  I_E''(\delta^{-1}),
\]
and by a standard computation and \eqref{E:logMEderivs},
\begin{equation} \label{E:IEderivs}
    I_E'(\delta^{-1}) = x_0, \qquad
    I_E''(\delta^{-1}) = \frac{1}{(\log M_E)''(x_0)} \sim \frac{|x_0|^{3-c}}
      {\Gamma(3-c)\varphi(|x_0|^{-1})} \quad \text{as } \delta \to 0.
      \end{equation}
It follows from these and \eqref{E:deltavsx0} that
\begin{equation} \label{E:Fsecond}
   \delta^2 F''(\delta) \sim K_{11}x_0(\delta).
   \end{equation}

Next observe that integrating \eqref{E:logMEderivs} with $k=1$ gives
\begin{equation} \label{E:happrox}
  \log M_E(-x) \sim -\frac{\Gamma(2-c)}{c-1} x^{c-1} \varphi(x^{-1}) \quad \text{as } -x \nearrow 0.
  \end{equation}
With \eqref{E:ma0} and \eqref{E:LDfacts} this shows that
\begin{equation} \label{E:betaDelta}
   \beta\Delta \sim K_{12} M^{-(c-1)}\varphi(M) \quad \text{as } \beta\Delta \to 0.
\end{equation}

Observe that by 
\eqref{E:happrox} and \eqref{E:deltavsx0} we have
\begin{equation} \label{E:IEasymp}
  I_E(\delta^{-1}) = \delta^{-1}x_0 - \log M_E(x_0)
      \sim \frac{\Gamma(3-c)}{c-1} |x_0|^{c-1} \varphi(|x_0|^{-1}) \quad \text{as } \delta \to 0,
      \end{equation}
so by \eqref{E:deltavsx0} and \eqref{E:x0vsdelta},
\begin{equation} \label{E:excurscost}
  \delta I_E(\delta^{-1}) \sim \frac{2-c}{c-1} |x_0|
    \sim K_{13} \delta^{1/(2-c)} \widehat{ \left( \frac{1}{\varphi} \right) }_{2-c}\left( \frac{1}{\delta} \right)
    \quad \text{as } \delta \to 0.
    \end{equation}
If we take $\delta = \delta^*$ in \eqref{E:excurscost}, we have $|x_0| = \alpha_0$ and thus
\begin{equation} \label{E:relsize}
  \beta\Delta\delta^* - \delta^* I_E((\delta^*)^{-1}) = \alpha_0, \quad
    \beta\Delta\delta^* \sim \frac{\alpha_0}{c-1}, \quad 
    \delta^* I_E((\delta^*)^{-1}) \sim \frac{2-c}{c-1} \alpha_0 \quad
    \text{as } \beta\Delta \to 0.
    \end{equation}
From \eqref{E:ma0} and \eqref{E:relsize} we have both
\begin{equation} \label{E:energylimit}
   \beta \Delta \delta^* M \to K_{14} \quad \text{as } \beta\Delta \to 0,
   \end{equation}
and
\begin{equation} \label{E:freeenbal}
   \beta f^a(\beta,u) = \alpha_0 \sim (c-1)
      \beta\Delta\delta^* \quad \text{as } \beta\Delta \to 0.
\end{equation}

From \eqref{E:betaDelta} we have 
\[
  \frac{1}{\beta\Delta} \sim K_{15} \frac{1}{\alpha_0^{c-1}} 
  \ophi_{c-1} \left( \frac{1}{\alpha_0^{c-1}} \right) \quad
  \text{as } \beta\Delta \to 0.
  \]
By \eqref{E:relsize} and Lemma 1.10 of \cite{Se76} this means that
\[
  \frac{1}{(\beta\Delta \delta^*)^{c-1}} \sim \frac{K_{16}}{\alpha_0^{c-1}} 
  \sim K_{17}\frac{1}{\beta\Delta} \ophi_{c-1}^*\left( \frac{1}{\beta\Delta} \right),
   \]
or equivalently
\[
  \delta^* \sim K_{18} (\beta\Delta)^{(2-c)/(c-1)}  \hat{\varphi}_{c-1}\left( \frac{1}{\beta\Delta} \right),
  \]
and then also, by \eqref{E:ma0} and \eqref{E:freeenbal},
\begin{equation} \label{E:alpha0val}
  \frac{K_{9}}{M} \sim
  \alpha_0 \sim K_{19} (\beta\Delta)^{1/(c-1)}  \hat{\varphi}_{c-1}\left( \frac{1}{\beta\Delta} \right).
  \end{equation}
Thus Theorem \ref{T:anneal} is proved.

For the remainder of this section we consider $c>3/2$.
Define $\Delta_1 = \Delta_1(\beta)$ to be the unique positive $\Delta$ where $\delta^*(\beta\Delta) = 2\Delta/\beta$ (i.e. where the linear upper bound intersects the contact fraction curve) and let $M_1 = M(\beta\Delta_1)$; see the discussion after Theorem \ref{T:bigc}.
Using \eqref{E:betaDelta} we get that for $\Delta \sim \Delta_1(\beta)$,
\[
  \frac{2\Delta}{\beta} \sim K_{12} \frac{\varphi(M)}{\beta^2 M^{c-1}}, \qquad 
    \delta^* \sim \delta_M \sim \frac{K_{20}}{M^{2-c}\varphi(M)} \quad \text{as } \beta \to 0,
  \]
so for $\Delta \sim \Delta_1$, or equivalently, $2\Delta/\beta \sim \delta^*$, we have 
\begin{equation} \label{E:Mbeta}
  \frac{1}{\beta^2} \sim K_{21} \frac{M^{2c-3}}{\varphi(M)^2} \quad \text{as } \beta \to 0,
  \end{equation}
and therefore from \eqref{E:betaDelta} again,
\begin{equation} \label{E:Delta1prop}
  \Delta_1^2 = \frac{(\beta\Delta_1)^2}{\beta^2} \sim \frac{K_{22}}{M_1}.
  \end{equation}
A similar argument shows that for $\Delta \ll \Delta_1$ we have $\Delta^2 \ll 1/M$, and for $\Delta \gg \Delta_1$ we have $\Delta^2 \gg 1/M$, confirming comments after Theorem \ref{T:bigc}.  Taking $\Delta = K\Delta_1$ for some $K$, we have from \eqref{E:LDfacts}, \eqref{E:alpha0val} and \eqref{E:Delta1prop} that 
\[
  \frac{\Delta^2}{2} \sim \frac{K^2 K_{22}}{2M(\beta\Delta_1)} 
    \sim \frac{K^{2-\frac{1}{c-1}} K_{22}}{2M(\beta\Delta)} \quad \text{as } \beta \to 0.
    \]
Letting $K_{23}$ be defined by $K_{23}^{2-\frac{1}{c-1}} K_{22} = 2K_{9}$, and letting
$\Delta_2 = K_{23} \Delta_1$, we see that $\Delta^2/2 \sim r\beta f^a(\beta,u)$ with 
$r<1$ if $K < K_{23}$, $r>1$ if $K > K_{23}$, and $r=1$ if $K=K_{23}$, i.e. if $\Delta \sim \Delta_2$.  Thus asymptotically, the annealed free energy curve intersects the upper bound $\Delta^2/2$ at approximately $\Delta = \Delta_2$.

To complete verification of the comments after Theorem \ref{T:bigc}, we show that $K_{3}$ can be chosen (in the definition of $\Delta_0$) so that $\Delta_0 \sim \Delta_2$.  From \eqref{E:alpha0val} and \eqref{E:Delta1prop} we have 
\[
  \frac{1}{\Delta_1^2} \sim \frac{K_{9}}{K_{22}K_{19}}
    \left( \frac{1}{\beta\Delta_1} \right)^{1/(c-1)} \ophi_{c-1}^* 
    \left( \frac{1}{\beta\Delta_1} \right)^{1/(c-1)} \quad \text{as } \beta \to 0,
    \]
or equivalently,
\[
  \frac{1}{\Delta_1^{2c-2}} \sim K_{24}
    \frac{1}{\beta\Delta_1} \ophi_{c-1}^* 
    \left( \frac{1}{\beta\Delta_1} \right),
    \]
or
\[
  \frac{1}{\Delta_1^{2c-2}} \ophi_{c-1}\left( \frac{1}{\Delta_1^{2c-2}} \right)
    \sim \frac{K_{24}}{\beta\Delta_1},
  \]
or from \eqref{E:Mbeta}  and \eqref{E:Delta1prop}, since $c>3/2$,
\[
  \frac{1}{\Delta_1^{2c-3}} \ophi_{c-\frac{3}{2}} \left( \frac{1}{\Delta_1^{2c-3}} \right) =
  \frac{1}{\Delta_1^{2c-3}} \varphi\left( \frac{1}{\Delta_1^2} \right)^{-1} \sim \frac{K_{24}}{\beta},
  \]
or
\[
  \frac{1}{\Delta_1^{2c-3}} \sim \frac{K_{24}}{\beta}
    \ophi_{c-\frac{3}{2}}^* \left( \frac{1}{\beta} \right) =
    \frac{K_{24}}{\beta} \hat{\varphi}_{c-\frac{3}{2}}  \left( \frac{1}{\beta} \right)^{-(c-\frac{3}{2})}.
  \]
Thus taking $K_{3} = K_{23}K_{24}^{-1/(2c-3)}$ we get 
\[
  \Delta_2 = K_{23} \Delta_1 \sim 
    K_{3} \beta^{1/(2c-3)} \hat{\varphi}_{c-\frac{3}{2}}  \left( \frac{1}{\beta} \right)^{1/2}
    = \Delta_0 \quad \text{as } \beta \to 0,
  \]
as desired.

\section{Proof of Theorem \ref{T:bigc}} \label{S:bigproof}
We may assume that $-\frac{\beta}{2} + \Delta \in \mD(\beta)$, i.e. a unique quenched contact fraction exists.

We begin with the proof that the quenched free energy is close to the annealed
when $\Delta$ is not too small, which is the first part of \eqref{E:highfree2}. 

For a disorder realization $\{v_j\}$ and Markov chain trajectories
$\{x_n\},\{y_n\}$ define
\[
         L_N(\{x_n\}) = \sum_{n=1}^N \delta_{\{ x_n = 0 \}}
\]
\[
         R_N(\{v_n\},\{x_n\}) = \sum_{n=1}^N (u + v_n)\delta_{\{x_n = 0\}}
\]
\[
         B_N(\{x_n\},\{y_n\}) = \sum_{n=1}^N \delta_{\{x_n = y_n = 0\}}.
\]
We write only $L_N, B_N$ when confusion is unlikely.

We first consider the free energy with a modified measure $\hP^X$ in place of $P^X$ in \eqref{E:GibbsWtQ}.  Let $M$ be the correlation length from \eqref{E:Mdef}.  $\hP^X$ is defined by the specification that under $\hP^X$, $E_1$ is uniform in $\{1,...,l_0M\}$ and $E_2,E_3,...$ are iid with distribution $P^X(E_1 \in \cdot)$, independent of $E_1$.  Here $l_0$ is an integer to be specified.  In a harmless abuse of notation, we denote the corresponding expectation and conditional expectation by $\langle \cdot \rangle^{\hat{X}}$ and $E^{\hat{X}}( \cdot | \cdot)$, respectively.  We use $\hZ_{[0,N]}^{\{V_i\}}(\beta,u)$ and $\hm_{[0,N]}^{\beta,u,\{V_i\}}$ to denote the corresponding quenched partition function and finite-volume Gibbs distribution, respectively. 
In the annealed case when $\hP^X$ is in effect we use the notation $\hZ_{[0,N]}(\beta,u)$ 
and $\hm^{\beta,u}_{[0,N]}$.

For $\Xi_N \subset \{0\} \times \Sigma^N$ a set of length-$N$ Markov chain 
trajectories with $x_0 = 0$, let
$Z^{\{V_i\}}_{[0,N]}(\beta,u,\Xi_N)$ denote the contribution to
$Z^{\{V_i\}}_{[0,N]}(\beta,u)$ from trajectories in $\Xi_N$, and similarly for $\hZ$ in place of $Z$.

Fix an integer $k_0 > 6l_0$ to be specified, and define a block length $N = k_0 M$.  The $k$th \emph{block} is then $((k-1)N,kN] \cap \ZZ$.  We would like to choose events $\Xi_{mN}, m \geq 1$ such that such that $Z^{\{V_i\}}_{[0,mN]}(\beta,u,\Xi_{mN})$ is a ``not too small'' fraction of $Z^{\{V_i\}}_{[0,mN]}(\beta,u)$, and $Z^{\{V_i\}}_{[0,mN]}(\beta,u,\Xi_{mN})$ approximately factors into a product of contributions from each block.  Such a product form can be made exactly true, as in \cite{AS06} and \cite{GT06}, by conditioning on a return to 0 at the end of every block, but the fact that the entropy cost of such conditioning grows large as $\Delta \to 0$ (i.e. as $N \to \infty$) makes it unworkable in our context.  Here we modify the product requirement slightly, to allow trajectories to bypass ``bad'' blocks.  To that end, in the $k$th block there is a \emph{landing zone} $((k-1)N,(k-1)N+l_0 M] \cap \ZZ$ at the beginning of the block, a \emph{prohibited zone} $(kN-l_0 M,kN] \cap \ZZ$ at the end of the block
and a \emph{takeoff zone} $(kN-5l_0 M,kN-l_0 M] \cap \ZZ$ just before the prohibited zone.  Let $m \geq 1$ and $J \subset \{1,..,m\}$ with $1 \in J$, and label the elements of $J$ as $j_1<..<j_{|J|}$.  For $m \geq 1$ we define $\Xi^J_{mN}$ to be the set of trajectories $x_{[0,mN]}$ satisfying the following criteria:
\begin{itemize}
\item[(a)] for each $2 \leq k \leq |J|$ there is an excursion which starts in the takeoff zone of block $j_{k-1}$ and ends in the landing zone of block $j_k$;
\item[(b)] for each $k \leq |J|$ there is at least one return to 0 in the first half of the takeoff zone of block $j_k$;
\item[(c)] there are no returns to 0 after the takeoff zone of block $j_{|J|}$.
\end{itemize}
We write $\Xi_N$ for $\Xi_N^{\{1\}}$.

Let $0 < \epsilon < 1/2$.  Suppose that we can choose $k_0,l_0$ such that for constants $K_i$ (not depending on $\epsilon$) to be specified, we have the following properties:
\begin{itemize}
\item[(P1)] The partition function approximately factors, in that for all $m \geq 1$ and $J \subset \{1,..,m\}$ with $1 \in J$,
\[
  \hZ_{[0,mN]}^{\{V_i\}}(\beta,u,\Xi^J_{mN}) \geq e^{-K_{25}l_0 m} \prod_{k \in J}
  \hZ_{[0,N]}^{\{V_{(k-1)N+i}\}}(\beta,u,\Xi_{N}).
  \]  
\item[(P2)]
\[
   \hm^{\beta,u}_{[0,N]} \left( \Xi_N \right) \geq e^{-K_{26}l_0}.
\]
\item[(P3)] For $\{X_n\}$ and $\{Y_n\}$ two
independent realizations of the Markov chain, 
\[
  \hm^{\beta,u}_{[0,N]} \left( e^{ \beta^2 B_N(\{X_i\},\{Y_i\}) } - 1\ \bigg|\  
    \{X_i\},\{Y_i\} \in \Xi_N \right) \leq \frac{\epsilon}{8}.
    \]
\end{itemize}
Here $\{V_{(k-1)N+i}\}$ denotes the sequence $\{ V_{(k-1)N+1},V_{(k-1)N+2},... \}$.
Note that (P3) says that two randomly
chosen trajectories from
$\Xi_N$ intersect each other at 0 only rarely.
Also, there is some competition between (P3) and (P2), in that for (P3) we would
like $N$ not too big so that $B_N$ is not too big, while for (P2) to be a workable lower bound, we must have the number $N/M$ of correlation lengths be much greater than $l_0$. 

Our basic technique for proving the first part of \eqref{E:highfree2} from (P1)--(P3) is the second moment method.  Note first that
\begin{equation} \label{E:mean}
  \left\langle \hZ_{[0,N]}^{\{V_i\}}(\beta,u,\Xi_N) \right\rangle^V
    = \left\langle e^{\beta\Delta L_N} \delta_{\Xi_N} \right\rangle^{\hat{X}}
    = \hm^{\beta,u}_{[0,N]}(\Xi_N) \left\langle e^{\beta\Delta L_N} \right\rangle^{\hat{X}}
    = \hm^{\beta,u}_{[0,N]}(\Xi_N) \hZ_{[0,N]}(\beta,u).
\end{equation}
We have for trajectories $\{x_i\}, \{y_i\}$ that
\begin{align}
  &\left\langle e^{\beta R_N(\{x_i\})} e^{\beta R_N(\{y_i\})} \right\rangle^V \notag \\
  &\qquad = \exp\big( \beta\Delta(L_N(\{x_i\}) - B_N(\{x_i\},\{y_i\})) \big)
    \exp\big( \beta\Delta(L_N(\{y_i\}) - B_N(\{x_i\},\{y_i\})) \big) \notag \\
  &\qquad \qquad \qquad \cdot \exp\big( (2\beta u + 2\beta^2)B_N(\{x_i\},\{y_i\}) \big) \notag \\
  &\qquad = e^{ \beta\Delta L_N(\{x_i\}) } e^{\beta\Delta L_N(\{y_i\}) } 
    e^{ \beta^2 B_N(\{x_i\},\{y_i\}) },
\end{align}
so by (P3), provided $\beta$ is small,
\begin{align} \label{E:variance}
  \var^V&\left( \hZ_{[0,N]}^{\{V_i\}}(\beta,u,\Xi_N) \right) \notag \\
  &= \sum_{\{x_n\} \in \Xi_N}\ \sum_{\{y_n\} \in \Xi_N}\ e^{ \beta\Delta L_N(\{x_i\}) } 
    e^{\beta\Delta L_N(\{y_i\}) } \left( e^{ \beta^2 B_N(\{x_i\},\{y_i\}) } - 1 \right) 
    \hP^X(\{x_i\}) \hP^X(\{y_i\}) \notag \\
  &= \left( \left\langle \hZ_{[0,N]}^{\{V_i\}}(\beta,u,\Xi_N) \right\rangle^V \right)^2
     \hm^{\beta,u}_{[0,N]} \left( e^{ \beta^2 B_N(\{X_i\},\{Y_i\}) } - 1\ \bigg|\  \{X_i\},\{Y_i\} \in \Xi_N    
     \right) \notag \\
   &\leq \frac{\epsilon}{8} \left( \left\langle \hZ_{[0,N]}^{\{V_i\}}(\beta,u,\Xi_N) \right\rangle^V \right)^2.
\end{align}
Also, by \eqref{E:mean} and (P2),
\begin{align}
  \left\langle \hZ_{[0,N]}^{\{V_i\}}(\beta,u,\Xi_N) \right\rangle^V
    &\geq K_{27} e^{-K_{26}l_0} \hZ_{[0,N]}(\beta,u).\notag
\end{align}
Therefore by Chebyshev's inequality and \eqref{E:variance},
\begin{equation} \label{E:Chebybound}
  P^V\left( \hZ_{[0,N]}^{\{V_i\}}(\beta,u,\Xi_N) \leq \half K_{27} e^{-K_{26}l_0} \hZ_{[0,N]}(\beta,u) \right)
    \leq \frac{\epsilon}{2}.
    \end{equation}
We say the $k$th block is \emph{good} if $k=1$ or
\[
  \hZ_{[0,N]}^{ \{V_{(k-1)N+i}\} }(\beta,u,\Xi_N) > \half K_{27} e^{-K_{26}l_0} \hZ_{[0,N]}(\beta,u),
  \]
and \emph{bad} otherwise, and we let $J_m = J_m(\{V_i\})= \{k \leq m:$ block $k$ is good$\}$.  By (P1),
\begin{align} \label{E:freelowbd}
  \frac{1}{mN} \log \hZ_{[0,mN]}^{\{V_i\}}(\beta,u) &\geq 
    \frac{1}{mN} \hZ_{[0,mN]}^{\{V_i\}}(\beta,u,\Xi^J_{mN}) \\
  &\geq \frac{\log(K_{27}/2)}{N} - (K_{26}+K_{25})\frac{l_0}{N} + \frac{|J_m|}{mN} \log 
     \hZ_{[0,N]}(\beta,u), \notag
\end{align}
so with $P^V$-probability one,
\begin{align} \label{E:freelowbd2}
 \beta f^q(\beta,u) &= \lim_{m \to \infty} \frac{1}{mN} \log \hZ_{[0,mN]}^{\{V_i\}}(\beta,u) \\
 &\geq \frac{\log(K_{27}/2)}{N}
   - (K_{26}+K_{25})\frac{l_0}{N} + \left( 1 - \frac{\epsilon}{2} \right) 
   \frac{1}{N} \log \hZ_{[0,N]}(\beta,u). \notag
\end{align}
To compare $N^{-1} \log \hZ_{[0,N]}(\beta,u)$ to $\beta f^a(\beta,u)$ we use subadditivity and the representation $\hZ_{[0,N]}(\beta,u) = \left\langle e^{\beta\Delta L_N} \right\rangle^{\hat{X}}$.  It is easily seen that for all $n,k$,
\begin{equation} \label{E:submult}
  \left\langle e^{\beta\Delta L_{n+k}} \right\rangle^X \leq \left\langle e^{\beta\Delta L_n} \right\rangle^X
    \left\langle e^{\beta\Delta (1+L_k)} \right\rangle^X,
    \end{equation}
so $a_n =  \beta\Delta + \log \left\langle e^{\beta\Delta L_n} \right\rangle^X, n \geq 1$, defines a subadditive sequence.  Therefore $a_k/k \geq \lim_n a_n/n = \beta f^a(\beta,u)$ for each fixed $k$.  Hence for each $u \leq l_0 M$ we have
\[
  E^{\hat{X}}( e^{\beta\Delta L_N} \mid E_1 = u) 
    = e^{\beta\Delta}  \left\langle e^{\beta\Delta L_{N-u}} \right\rangle^X
    \geq e^{(N-l_0 M)\beta f^a(\beta,u)},
  \]
and therefore
\[
  \hZ_{[0,N]}(\beta,u) = \left\langle e^{\beta\Delta L_N} \right\rangle^{\hat{X}} 
    \geq e^{(N-l_0 M)\beta f^a(\beta,u)}.
  \]
It follows that
\begin{equation} \label{E:Nfreebd}
  \frac{1}{N} \log \hZ_{[0,N]}(\beta,u) \geq \left( 1 - \frac{l_0}{k_0} \right) \beta f^a(\beta,u).
  \end{equation}
From \eqref{E:LDfacts} and \eqref{E:ma0} we have 
\[
  \fN \leq \frac{l_0}{N} \leq K_{29}\frac{l_0}{k_0}\beta f^a(\beta,u),
  \]
which with \eqref{E:freelowbd2} and \eqref{E:Nfreebd} shows that, provided $k_0/l_0$ is sufficiently large (depending on $\epsilon$), we have
\[
  \beta f^q(\beta,u) \geq (1-\epsilon)\beta f^a(\beta,u),
  \]
completing the proof of the first part of \eqref{E:highfree2}.

Our remaining task to prove the first part of \eqref{E:highfree2}, then, is to establish (P1)--(P3).  We begin with (P1).  Fix $m,J$ as is (P1), and assume $|J| \geq 2$.  Let $U_k$ and $T_k$ denote the location of the first and last returns, respectively, in the $k$th block (when such exist), necessarily in the landing and takeoff zones, respectively, for trajectories in $\Xi^J_{mN}$ and good blocks $k$.  Suppose $j_{|J|} < m$; we then first deal with the final excursion initiated before time $mN$, necessarily from the takeoff zone of the last good block, as follows.  We have
\begin{align} \label{E:lastexc}
  \hZ_{[0,mN]}^{\{V_i\}}(\beta,u,\Xi^J_{mN}) &= \sum_t 
    \hZ_{[0,j_{|J|}N]}^{\{V_i\}}(\beta,u,\Xi^J_{j_{|J|}N} \cap \{ T_{j_{|J|}} = t \} ) 
      \frac{ P^X(E_1 > mN-t) }{ P^X(E_1 > j_{|J|}N - t) }, 
\end{align}
where the sum is over $t$ in the takeoff zone of block $j_{|J|}$.  Since $j_{|J|}N - t \geq l_0 M$ and $mN - t \leq (m - j_{|J|})k_0 M + 5l_0 M$, there exist constants $K_i$ such that provided $l_0 \geq K_{30}$,
$k_0 \leq l_0 e^{K_{31}l_0}$ and $\beta\Delta$ is sufficiently small 
(depending on $\epsilon$) so that $M$ is large, we have
\[
  \frac{ P^X(E_1 > mN-t) }{ P^X(E_1 > j_{|J|}N - t) } \geq \half
    \left( \frac{l_0}{(m - j_{|J|})k_0 + 5l_0} \right)^{c-1}
    \geq e^{-K_{32}(m - j_{|J|})l_0}.
  \]
Then by \eqref{E:lastexc},
\begin{equation} \label{E:lastexc2}
  \hZ_{[0,mN]}^{\{V_i\}}(\beta,u,\Xi^J_{mN}) \geq e^{-K_{32}(m - j_{|J|})l_0}
    \hZ_{[0,j_{|J|}N]}^{\{V_i\}}( \beta,u,\Xi^J_{j_{|J|}N} ).
    \end{equation}
Having effectively replaced $m$ with $j_{|J|}$,
we next decompose according to the starting and ending points $t,u$ of the excursion from block $j_{|J|-1}$ to block $j_{|J|}$, as follows:
\begin{align} \label{E:Zdecomp}
  \hZ_{[0,j_{|J|}N]}^{\{V_i\}}(\beta,u,\Xi^J_{j_{|J|}N}) &= \sum_{t,u} 
    \hZ_{[0,j_{|J|-1}N]}^{\{V_i\}}(\beta,u,\Xi^J_{j_{|J|-1}N} \cap \{ T_{j_{|J|-1}} = t \} ) \\
  &\qquad \qquad \cdot \hZ_{[0,N]}^{ \{ V_{(j_{|J|} - 1)N + i} \} } 
    (\beta,u,\Xi_N \cap \{ U_1 = u - (j_{|J|} - 1)N \} ) \notag \\
  &\qquad \qquad \cdot \frac{ P^X(E_1 = u-t) }{ P^X( E_1 > j_{|J| - 1}N - t) \hP^X( E_1 = u - (j_{|J|} - 1)N ) }.
    \notag 
\end{align}
Note that
\[
  \left( l_0 + (j_{|J|} - j_{|J|-1} - 1)k_0 \right) M \leq u-t \leq \left( 6l_0 + (j_{|J|} - j_{|J|-1} - 1)k_0 \right) M, 
  \]
$j_{|J|-1}N - t \geq l_0 M$ and $\hP^X( E_1 = u - (j_{|J|} - 1)N ) = 1/l_0 M$.  Assuming again that $l_0 \geq K_{30}$, $k_0 \leq l_0 e^{K_{31}l_0}$ and $\beta\Delta$ is sufficiently small, it follows readily that
\begin{align} \label{E:fracbound}
  \frac{ P^X(E_1 = u-t) }{ P^X( E_1 > j_{|J| - 1}N - t) \hP^X( E_1 = u - (j_{|J|} - 1)N ) } &\geq K_{33}
    \left( 6 + ( j_{|J|} - j_{|J|-1} - 1) \frac{k_0}{l_0} \right)^{-c} \notag \\
  &\geq e^{-K_{34}( j_{|J|} - j_{|J|-1} )l_0}.
  \end{align}
Note it is essential that the correlation length $M$ cancels out in \eqref{E:fracbound}, so that the lower bound does not depend on $\beta\Delta$.  From \eqref{E:Zdecomp} and \eqref{E:fracbound},
\[
  \hZ_{[0,j_{|J|}N]}^{\{V_i\}}(\beta,u,\Xi^J_{j_{|J|}N}) \geq e^{-K_{34}( j_{|J|} - j_{|J|-1} )l_0} 
    \hZ_{[0,j_{|J|-1}N]}^{\{V_i\}}(\beta,u,\Xi^J_{j_{|J|-1}N} ) \hZ_{[0,N]}^{ \{ V_{(j_{|J|} - 1)N + i} \} } 
    (\beta,u,\Xi_N).
    \]
Iterating this we obtain
\[
  \hZ_{[0,j_{|J|}N]}^{\{V_i\}}(\beta,u,\Xi^J_{j_{|J|}N}) \geq e^{-K_{34}j_{|J|} l_0} 
    \prod_{k \in J} \hZ_{[0,N]}^{\{V_{(k-1)N+i}\}}(\beta,u,\Xi_{N}),
    \]
which with \eqref{E:lastexc2} completes the proof of (P1).

We next prove (P2), for which we need only consider one length-$N$ block.  Let $D_N$ denote the event that there is at least one return to 0 in the first half of the takeoff zone, i.e. in $(N - 5l_0 M,N - 3l_0 M]$, and let $C_N$ denote the event that there are no returns to 0 after the takeoff zone, i.e. after $N - l_0 M$, so $\Xi_N = D_N \cap C_N$.  We introduce the tilted measure $Q=Q^X$, on trajectories of the Markov chain, given by
\begin{align} \label{E:QXdef}
   Q^X&(E_1 = k_1,..,E_m = k_m) \\
   &= \frac{e^{-\alpha_0(k_1 + ..+ k_m)}}{(\langle
     e^{-\alpha_0 E_1} \rangle^X)^m} P^X(E_1 = k_1,..,E_m = k_m) \quad
     \text{ for all } m,k_1,..,k_m, \notag
\end{align}
which by \eqref{E:alphadef} satisfies
\begin{equation} \label{E:Qmean}
   \langle E_1 \rangle^Q = \frac{1}{\delta^*},
\end{equation}
where $\langle \cdot \rangle^Q$ denotes expected value under $Q^X$.  Equation
\eqref{E:QXdef} does not of course completely determine a distribution for
trajectories--we complete the definition by specifying that the
conditional distribution of the trajectory under $Q^X$ given $E_1 = k_1, E_2 =
k_2,...$ is the same as the corresponding conditional distribution under $P^X$.  Observe that $\hm^{\beta,u}_{[0,N]}$ tilts $P^X$ so that the contact fraction becomes $\delta^*$, i.e. so that the mean excursion length is approximately $1/\delta^*$ for large $N$, so from \eqref{E:Qmean} we expect the measures $\hm^{\beta,u}_{[0,N]}$ and $Q^X$ to be similar, a relation we will now make precise.  Given $n \leq k \leq r$ we have from \eqref{E:LDfacts} that
\begin{equation} \label{E:measures}
  \mu^{\beta,u}_{[0,r]}(\tau_n = k) = \frac{ e^{\beta\Delta n} \left\langle e^{\beta\Delta L_{r-k}} \right\rangle^X }
    {  \left\langle e^{\beta\Delta L_r} \right\rangle^X } P^X(\tau_n = k)
    = \frac{ e^{\alpha_0 k} \left\langle e^{\beta\Delta L_{r-k}} \right\rangle^X }
    {  \left\langle e^{\beta\Delta L_r} \right\rangle^X } Q^X(\tau_n = k).
    \end{equation}
Fix $u \leq l_0 M$ and let $v = N - 4l_0M - u$.  We now specify $n = \lfloor v\delta^* \rfloor$.  Then by \eqref{E:measures},
\begin{align} \label{E:DNlower}
  \hm^{\beta,u}_{[0,N]}(D_N \mid E_1 = u) &\geq \mu^{\beta,u}_{[0,N-u]}\bigg( 
    (k_0 - 5l_0)M - u < \tau_n \leq (k_0 - 3l_0)M - u \bigg) \notag \\
  &\geq \frac{ e^{\alpha_0 (N - 6l_0 M)} }{ \left\langle e^{\beta\Delta L_N} \right\rangle^X } Q^X\bigg( 
    v - l_0 M < \tau_n \leq v + l_0 M \bigg). 
  \end{align}
By \eqref{E:Qmean} we have
\[
  \langle \tau_n \rangle^Q \sim v \quad \text{as } \beta\Delta \to 0,
  \]
and the variance of $E_1$ under $Q^X$ is easily shown (using \eqref{E:ma0}) to satisfy
\begin{equation} \label{E:varE1}
   \var^Q(E_1) \sim \langle E_1^2 \rangle^Q 
   = e^{\beta\Delta} \sum_{j=1}^\infty e^{-\alpha_0 j} j^{2-c} \varphi(j)
   \sim K_{35}M^{3-c}\varphi(M),
\end{equation}
so using \eqref{E:deriv1},
\begin{equation} \label{E:vartauJ}
   \var^Q(\tau_n) \sim K_{35}\delta^* v M^{3-c}\varphi(M) \leq K_{36} k_0 M^2,
\end{equation}
also as $\beta\Delta \to 0$.  
Therefore provided $k_0 \leq l_0^2/4K_{36}$, by Chebyshev's
inequality and \eqref{E:DNlower} we have for $\beta\Delta$ small that
\begin{equation} \label{E:DNlower2}
  \hm^{\beta,u}_{[0,N]}(D_N) \geq \frac{ e^{\alpha_0 (N - 6l_0 M)} }
    { 2\left\langle e^{\beta\Delta L_N} \right\rangle^X }.
    \end{equation}

Let $G_N$ denote the time of the first return to 0 in the takeoff zone (necessarily in $(N - 5l_0 M,N - 3l_0 M]$, if the event $D_N$ occurs), when such a return exists.  For $g \in (N - 5l_0 M,N - 3l_0 M]$ we have using \eqref{E:deltan} that 
\begin{align} \label{E:CNcondl}
  \hm^{\beta,u}_{[0,N]}&(C_N \mid G_N = g) \notag \\
  &\geq \mu^{\beta,u}_{[0,N-g]}(\tau_i = k, L_{N-g} = i \text{ for some } k \leq N-g-l_0 M
    \text{ and } i \geq 1) \\
  &= \sum_{k \leq N-g-l_0 M} \sum_{i \geq 1} \frac{ e^{\beta\Delta i} }
    { \left\langle e^{\beta\Delta L_{N-g}} \right\rangle^X }
    P^X(\tau_i = k) P^X(E_1 > N-g-k) \notag \\
  &\geq \frac{ P^X(E_1 > 5l_0 M) }{ \left\langle e^{\beta\Delta L_{5l_0 M}} \right\rangle^X } 
    \sum_{k \leq 2l_0 M} P^X(X_k = 0) \notag \\
  &\geq \frac{P^X(E_1 > 5l_0 M) \left\langle L_{2l_0 M} \right\rangle^X}
    { \left\langle e^{\beta\Delta L_{5l_0 M}} \right\rangle^X } \notag \\
  &\geq \frac{ K_{37}}
    { \left\langle e^{\beta\Delta L_{5l_0 M}} \right\rangle^X }. \notag 
  \end{align}
Combining this with \eqref{E:ma0} and \eqref{E:DNlower2} we get
\begin{align} \label{E:CNDN}
  \hm^{\beta,u}_{[0,N]}(\Xi_N) &= \hm^{\beta,u}_{[0,N]}(D_N) \hm^{\beta,u}_{[0,N]}(C_N \mid D_N) \\
  &\geq \frac{ K_{37} e^{\alpha_0 (N - 6l_0 M)} }
    {  2\left\langle e^{\beta\Delta L_N} \right\rangle^X 
    \left\langle e^{\beta\Delta L_{5l_0 M}} \right\rangle^X }. \notag
  \end{align}
  
We claim that for some $K_{38}$,
\begin{equation} \label{E:claim1}
  \left\langle e^{\beta\Delta L_{kM}} \right\rangle^X \leq K_{38} ke^{\alpha_0 kM} \quad \text{ for all } 
  k \geq 1.
  \end{equation}
Presuming this is proved, we have from \eqref{E:CNDN}, presuming once more that $k_0 \leq l_0 e^{K_{31}l_0}$ and $l_0$ is sufficiently large,
\begin{equation} \label{E:DNlower3}
  \hm^{\beta,u}_{[0,N]}(\Xi_N) \geq \frac{K_{37}}{50 K_{38}^2 k_0 l_0^2} 
    \exp\bigg( \alpha_0 (N - 6l_0 M) - \alpha_0 N -10 l_0 \alpha_0 M \bigg)
    \geq e^{-K_{39}l_0},
    \end{equation}
so (P2) is proved.  It should be pointed out that the various conditions we have required on $l_0$ and $k_0$ are compatible and can be summarized as follows:  there exist $K_{40}, K_{41}, K_{42}$ and $K_{43} = K_{43}(\epsilon)>1$ such that $k_0,l_0$ must satisfy
\[
  K_{40} \leq K_{43}l_0 \leq k_0 \leq K_{41}\min( l_0^2,l_0 e^{K_{42}l_0}).
  \]
To prove \eqref{E:claim1}, fix $k \geq 1$ and observe that $F(\delta) = \beta\Delta \delta - \delta I_E(\delta^{-1})$ is concave with $F(0) = 0$ and maximum value $F(\delta^*) = \alpha_0>0$, so $\delta_0$ given by $I_E(\delta_0^{-1}) = \beta\Delta$ is the unique positive solution of $F(\delta) = 0$.  By \eqref{E:IEderivs} we have $F'(\delta_0) =\delta_0^{-1}x_0(\delta_0)$.  By concavity $F$ is below its tangent line at $\delta_0$, so we have the bound
\begin{equation} \label{E:Fbound}
  F(\delta) \leq \begin{cases} \alpha_0 = |x_0(\delta^*)|, &\delta \leq \delta_0,\\ 
    \frac{x_0(\delta_0)}{\delta_0}(\delta - \delta_0), &\delta>\delta_0. \end{cases}
    \end{equation}
Let $j_0 = \lfloor \delta_0 kM \rfloor$.  Using the equivalence of $L_{kM} \geq j$ and $\tau_j \leq kM$ we obtain from \eqref{E:Fbound} that
\begin{align} \label{E:weightbd}
  \left\langle e^{\beta\Delta L_{kM}} \right\rangle^X &= 
    1 + \beta\Delta \sum_{j=1}^{kM} e^{\beta\Delta j} P^X(\tau_j \leq kM) \notag \\
  &\leq 1 + \beta\Delta \sum_{j=1}^{kM} e^{F(j/kM)kM} \\
  &\leq 1 + \beta\Delta j_0 e^{\alpha_0 kM} + 
    \beta\Delta \sum_{j >j_0} \exp\left( \frac{ x_0(\delta_0) }{ \delta_0 }
    \left( \frac{j}{kM} - \delta_0 \right) kM \right) \notag \\
  &\leq 1 + \beta\Delta\delta_0 kM e^{\alpha_0 kM} + 
    \beta\Delta \sum_{r=2}^\infty \sum_{(r-1)\delta_0 kM < j \leq r\delta_0 kM}
    e^{(r-2)x_0(\delta_0)kM} \notag \\
  &\leq 1 + \beta\Delta\delta_0 kM (e^{\alpha_0 kM} + K_{44}), \notag
  \end{align}
where the last inequality follows from $x_0(\delta_0) \leq -\alpha_0$ together with \eqref{E:ma0}.  Now we need to compare $\delta_0$ to $\delta^*$.  From \eqref{E:excurscost} and \eqref{E:relsize} we have
\begin{align} 
  \frac{\delta_0}{\delta^*} &= \frac{\beta\Delta\delta_0}{\beta\Delta\delta^*} \notag \\
  &\sim (2-c)\frac{x_0(\delta_0)}{x_0(\delta^*)} \notag \\
  &\sim (2-c)\left( \frac{\delta_0}{\delta^*} \right)^{1/(2-c)} \frac{ \psi(\delta_0^{-1}) }
    { \psi((\delta^*)^{-1}) }, \notag
\end{align}
with $\psi$ slowly varying, which implies that 
\[
   \frac{\delta_0}{\delta^*} \sim (2-c)^{-(2-c)/(c-1)}.
   \]
With \eqref{E:weightbd} and \eqref{E:energylimit} this proves the claim \eqref{E:claim1}, completing the proof of (P2).

Next we prove (P3).  We use $\langle \cdot \rangle^\mu$ to denote expectation under a measure $\mu$.  By (P2),
\begin{align} \label{E:removecond}
   \hm^{\beta,u}_{[0,N]} &\left( e^{ \beta^2 B_N(\{X_i\},\{Y_i\}) } - 1\ \bigg|\  
    \{X_i\},\{Y_i\} \in \Xi_N \right) \\
   &\leq \frac{1}{ \hm^{\beta,u}_{[0,N]}(\Xi_N)^2 } 
     \left\langle e^{ \beta^2 B_N(\{X_i\},\{Y_i\}) } - 1 \right\rangle^{\hm^{\beta,u}_{[0,N]}} \notag \\
   &\leq  e^{2K_{26} l_0} 
     \left\langle e^{ \beta^2 (B_N(\{X_i\},\{Y_i\}) + 1) } - 1 \right\rangle^{\mu^{\beta,u}_{[0,N]}}. \notag
   \end{align}
We can shift the problem from length scale $N$ to the length scale $M$, on which the measures $\mu^{\beta,u}_{[0,M]}$ and $P^X$ are comparable, via the inequality
\begin{equation} \label{E:MNcompare}
  \left\langle e^{ \beta^2 B_N(\{X_i\},\{Y_i\}) }\right\rangle^{\mu^{\beta,u}_{[0,N]}} \leq
    \left( \left\langle e^{ \beta^2 (B_M(\{X_i\},\{Y_i\}) + 1) } 
    \right\rangle^{\mu^{\beta,u}_{[0,M]}} \right)^{k_0}.
  \end{equation}
To quantify the comparability, by \eqref{E:claim1} (with trivial modifications to deal with the 2 in the exponent) and \eqref{E:ma0}, using the fact that $(x-1)^2 \leq x^2 - 1$ for $x \geq 1$ we have
\begin{align} \label{E:Holder}
   &\left( \left\langle e^{ \beta^2 (B_M(\{X_n\},\{Y_n\}) + 1) } - 1
     \right\rangle^{\mu^{\beta,u}_{[0,M]}} \right)^2 \\
   &\qquad \leq \left( \left\langle e^{\beta\Delta L_M(\{X_n\})}
     e^{\beta\Delta L_M(\{Y_n\})} \left( e^{ \beta^2 (B_M(\{X_n\},\{Y_n\}) + 1) } - 1
     \right) \right\rangle^X \right)^2 \notag \\
   &\qquad \leq \left\langle e^{2\beta\Delta L_M(\{X_n\})}
     \right\rangle^X \left\langle
     e^{2\beta\Delta L_M(\{Y_n\})} \right\rangle^X \left\langle
     e^{ 2\beta^2 (B_M(\{X_n\},\{Y_n\}) + 1) } - 1 \right\rangle^X \notag \\
   &\qquad \leq K_{45} \left\langle e^{ 2\beta^2 (B_M(\{X_n\},\{Y_n\}) + 1) } - 1 \right\rangle^X. \notag
\end{align}
Combining \eqref{E:MNcompare} and \eqref{E:Holder} we obtain
\begin{align} \label{E:NMcompare}
  &\left\langle e^{ \beta^2 (B_N(\{X_i\},\{Y_i\}) + 1) } - 1 \right\rangle^{\mu^{\beta,u}_{[0,N]}} \\
  &\qquad \leq e^{\beta^2} \left( \left\langle e^{ \beta^2 (B_M(\{X_i\},\{Y_i\}) + 1) } 
    \right\rangle^{\mu^{\beta,u}_{[0,M]}} \right)^{k_0} - 1 \notag \\
  &\qquad \leq e^{\beta^2} \left( 1 + K_{45}^{1/2} 
     \left( \left\langle e^{ 2\beta^2 (B_M(\{X_n\},\{Y_n\}) + 1) } - 1
     \right\rangle^X \right)^{1/2} \right)^{k_0} - 1. \notag 
\end{align}
To bound the expectation on the right side of \eqref{E:NMcompare} we will need some lemmas.
The first is a result of Garsia and Lamperti \cite{GL63} on renewal processes, which we specialize here to our situation.

\begin{lemma} \label{L:Erickson} (\cite{GL63})
Suppose that for some slowly varying function $\varphi$
and constants $K < \infty,\ 1<c<2$, the excursion length distribution of
an aperiodic Markov chain $\{X_n\}$ satisfies
\[
      P(E_1 = n) \leq K n^{-c}\varphi(n) \quad \text{ for all } n
\]
and
\[
      P(E_1 > n) \sim n^{-(c-1)}\varphi(n) \quad \text{as } n \to \infty.
\]
Then
\[
      P(X_n = 0) \sim \frac{\Gamma(2-c)}{\Gamma(c-1)} n^{-(2-c)}\varphi(n)^{-1}
         \quad \text{as } n \to \infty.
\]
\end{lemma}

Let $\tE_i$ denote the length of the $i$th excursion from $(0,0)$, $\tL_n$ the number of returns to 0
by time $n$ and $\tS_n$ the time of the $n$th return to $(0,0)$, for the chain $\{(X_i,Y_i)\}$.

\begin{lemma} \label{L:doubletail}
Let $\{X_i\}$ and $\{Y_i\}$ be independent copies of an aperiodic
recurrent Markov
chain starting at 0 and satisfying \eqref{E:tails} with $3/2 < c < 2$, or with $c=3/2$ and $\sum_{i=1}^\infty i^{-1} \varphi(i)^{-2}$ $= \infty$.  Then there exist
constants $K_i$ such that for all sufficiently large $n$,
\begin{equation} \label{E:doubletail}
    P(\tE_1 > n) \geq \frac{K_{46}}{\left\langle \tL_n \right\rangle^X} 
      \geq \begin{cases} K_{47} n^{-(2c-3)} \varphi(n)^2 &\text{if } c > \frac{3}{2},\\
      \frac{K_{48}}{ \sum_{i=1}^n i^{-1}\varphi(i)^{-2} } &\text{if } c=\frac{3}{2}. \end{cases}
    \end{equation}
\end{lemma}
\begin{proof}
We use $E_1, E_2,...$ to denote excursion lengths for $\{X_i\}$, as usual.  Define
\[
  \sigma_1 = \min\{i \geq 1: E_i > n\}, \qquad \sigma_2 = \min\{i \geq 1: \tE_i > n\}.
  \]
Let $T_{n1}$ and $T_{n2}$ denote the starting times of excursions $\sigma_1$ and $\sigma_2$ for the chains $\{X_i\}$ and  $\{(X_i,Y_i)\}$, respectively.  Then
\begin{equation} \label{E:sigmabound}
  \left\langle \tL_n \right\rangle^X \geq \left\langle \tL_n \delta_{\{T_{n2} \leq n/2\}} \right\rangle^X
    = \left\langle \sigma_{2} \delta_{\{T_{n2} \leq n/2\}} \right\rangle^X.
    \end{equation}
Suppose we can show that for some $K_{49}>0$,
\begin{equation} \label{E:claim2}
  P^X\left(T_{n2} \leq \frac{n}{2}\right) \geq K_{49} \quad \text{for all sufficiently large } n.
  \end{equation}	
Since $\sigma_2$ has a geometric distribution, this shows that 
\[
  \left\langle \sigma_{2} \delta_{\{T_{n2} \leq n/2\}} \right\rangle^X \geq K_{50}
    \left\langle \sigma_{2} \right\rangle^X = \frac{K_{50}}{P^X(\tE_1 > n)},
    \]
which with \eqref{E:sigmabound} completes the proof of the first inequality in \eqref{E:doubletail}.  The second inequality is a consequence of Lemma \ref{L:Erickson} and the relation
\[
   \left\langle \tL_n \right\rangle^X = \sum_{i=1}^n P^X(X_i = 0)^2.
   \]
To prove \eqref{E:claim2} we observe that
\[
  P^X\left(T_{n2} \leq \frac{n}{2}\right) \geq P^X\left(T_{n1} \leq \frac{n}{2}\right),
  \]
and proceed analogously to \eqref{E:CNcondl}.  Using \eqref{E:deltan} we have 
\begin{align} \label{E:Tn1bound}
  P^X\left(T_{n1} \leq \frac{n}{2}\right) &= \sum_{j \geq 0} P^X\left(\tau_j \leq \frac{n}{2}\right)P^X(E_{j+1}>n) \\
  &\geq P^X(E_1 > n) \left\langle L_{n/2} \right\rangle^X \notag \\
  &\geq K_{51}, \notag
  \end{align}
so \eqref{E:claim2} is proved.
\end{proof}

\begin{lemma} \label{L:boundreturns}
Let $\{X_i\}$ and $\{Y_i\}$ be independent copies of an aperiodic
recurrent Markov
chain starting at 0 and satisfying \eqref{E:tails} with $1 < c < 2$.  Then 
\begin{itemize}
\item[(i)] if $c > 3/2$ then there exists $K_{52} < \infty$ such that
\begin{equation}\label{E:bigc}
     P(B_N \geq k) \leq \left( 1 -
       \frac{\varphi(N)^2}{K_{52} N^{2c-3}} \right)^k
       \quad \text{for all } N,k \geq 1;
\end{equation}
\item[(ii)] if $c = 3/2$ and $\sum_n n^{-1} \varphi(n)^{-2}
= \infty$, then there exists $K_{53} < \infty$ such that
\begin{equation} \label{E:c32}
     P(B_N \geq k) \leq \left( 1 -
       \frac{1}{K_{53}\tphi(N)} \right)^k \quad \text{for all } N,k \geq 1;
\end{equation}
\item[(iii)] if $c < 3/2$, or if $c=3/2$ and $\sum_n n^{-1} \varphi(n)^{-2}
< \infty$,
then there exists $\epsilon_1 > 0$ such that
\begin{equation} \label{E:smallc}
     P(B_N \geq k) \leq \left( 1 -
       \epsilon_1 \right)^k \quad \text{for all } N,k \geq 1.
\end{equation}
\end{itemize}
\end{lemma}

\begin{proof}
By Lemma
\ref{L:Erickson} if $1<c<2$ we
have for some $K_{54}$
\begin{equation} \label{E:doublertn}
     P\big( (X_i,Y_i) = (0,0) \big) \leq K_{54} i^{-(4-2c)}\varphi(i)^{-2} \quad
      \text{for all } i \geq 1.
\end{equation}

If $1<c < 3/2$, or if $c=3/2$ and $\sum_n n^{-1} \varphi(n)^{-2} < \infty$,
then by \eqref{E:doublertn} we have $\sum_i\ P\big( (X_i,Y_i) = (0,0) \big) < \infty$ so the chain
$\{(X_i,Y_i)\}$ is transient and \eqref{E:smallc} follows.

If $c > 3/2$, then by Lemma \ref{L:doubletail}, we have
\[
    P(B_N > k) \leq P\left( \max_{j \leq k} \tE_j \leq N \right)
      \leq \left( 1 - \frac{\varphi(N)^2}{K_{52} N^{2c-3}} \right)^k,
\]
which proves \eqref{E:bigc}.
Similarly, if $c=3/2$ and $\sum_n n^{-1} \varphi(n)^{-2} = \infty$ then Lemma \ref{L:doubletail} gives \eqref{E:c32}.  
\end{proof}

We can now continue with the bound on the expectation on the right side of \eqref{E:NMcompare}.
By Lemma \ref{L:boundreturns}(i) we have
$B_M(\{X_n\},\{Y_n\}) + 1$ stochastically smaller (under $P^X$) than 
a geometric
random variable with parameter of form $p_M = 
K_{52}^{-1}M^{-(2c-3)}\varphi(M)^2$.  Let
\[
  a = \frac{ \epsilon e^{-2K_{26}l_0} }{ 32k_0 K_{45}^{1/2} },
  \]
where $K_{26}$ is from (P2) and $K_{45}$ from \eqref{E:Holder}, and suppose that, for $K_{55}$ to be specified,
\begin{equation} \label{E:Deltalower}
  \Delta \geq K_{55}a^{-(2c-2)/(2c-3)} \beta^{1/(2c-3)} 
    \hat{\varphi}_{c-\frac{3}{2}}\left( \frac{1}{\beta} \right)^{1/2},
  \end{equation}
which is a version of the assumption in \eqref{E:highfree2}.  From \eqref{E:Mbeta} and the discussion following it, for each fixed $K$ there exists $g(K)$, with $g(K) \nearrow \infty$ as $K\to \infty$, such that if we let $\Delta = \Delta(\beta) \to 0$, then the statement that $\Delta \sim K \Delta_0$ is equivalent to the statement that $p_M \sim g(K)\beta^2$ as $\beta \to 0$.  Thus if $K_{55}$ is large enough then
\begin{equation} \label{E:pMlower}
  p_M \geq \frac{4\beta^2}{a^2} \geq \left( 1 + \frac{1}{a^2} \right) \left( 1 - e^{-2\beta^2} \right),
  \end{equation}
so from the bound by a geometric random variable,
\begin{align}
    \left\langle e^{ 2\beta^2 (B_M(\{X_n\},\{Y_n\})+1) }
     - 1 \right\rangle^X &\leq \frac{e^{2\beta^2} - 1}{1 - (1-p_M)e^{2\beta^2}} \leq a^2.
\end{align}
Plugging this into \eqref{E:NMcompare} we obtain from \eqref{E:removecond} that provided $\beta$ is small enough (depending on $l_0$),
\begin{align}
  \hm^{\beta,u}_{[0,N]} \left( e^{ \beta^2 B_N(\{X_i\},\{Y_i\}) } - 1\ \bigg|\  
    \{X_i\},\{Y_i\} \in \Xi_N \right) &\leq e^{2K_{26}l_0} \left[ \exp\left( \beta^2 + 
    \frac{\epsilon e^{-2K_{26}l_0}}{32} \right) - 1 \right] \notag \\
  &< \frac{\epsilon}{8}.
\end{align}
The proof of (P3), and thus of the free-energy inequality in \eqref{E:highfree2}, is now complete.

We next consider to the contact-fraction inequality in \eqref{E:highfree2}.
Recall that for $F$ from \eqref{E:Fdef}, we have $F$ maximized at $\delta^*$ with 
$x_0(\delta^*) = -\alpha_0 = -F(\delta^*)$ (see \eqref{E:LDfacts}.)
Hence from \eqref{E:x0vsdelta} and \eqref{E:Fsecond}, for all $|\gamma| < \half$, we have as $\beta\Delta \to 0$ that
\[
   -(1+\gamma)^2 (\delta^*)^2 F''((1+\gamma)\delta^*) \sim
     (1+\gamma)^{1/(2-c)} K_{56}F(\delta^*).
\]
This shows that $-\delta^2 F''(\delta)/F(\delta^*)$ is bounded away
from 0 on $[\delta^*/2,3\delta^*/2]$,
uniformly for small $\beta\Delta$.  It follows that given $0<\lambda < 1/2$
there exists $\theta>0$ such that
\begin{equation} \label{E:Fconcavity}
   |\delta - \delta^*| \geq \lambda \delta^* \quad \text{implies }
   \quad F(\delta) \leq (1-\theta)F(\delta^*).
   \end{equation}
Fix $0<\lambda < \half$ and define the events
\[
   \Phi_{1N} = \left\{\{x_n\}: \frac{L_N(\{x_n\})}{N} \leq
     (1-\lambda)\delta^* \right\}, \quad
     \Phi_{2N} = \left\{\{x_n\}: \frac{L_N(\{x_n\})}{N} \geq
     (1+\lambda)\delta^* \right\}.
   \]
 From trivial modifications of Theorem 2.1 of \cite{AS06} and from \eqref{E:Fconcavity},
 we obtain that the
contributions to the
quenched and annealed free energy from trajectories in $\Phi_{iN}$ satisfy
\begin{align} \label{E:annealmax}
    \limsup_N &\fN \left\langle \log Z_{[0,N]}^{\{V_i\}}(\beta,u,\Phi_{iN})
     \right\rangle^V \\
   &\leq \limsup_N \fN \log \left\langle Z_{[0,N]}^{\{V_i\}}(\beta,u,\Phi_{iN})
     \right\rangle^V \notag \\
   &\leq \sup\{F(\delta):
     |\delta - \delta^*| > \lambda \delta^*\} \notag \\
   &\leq (1-\theta)F(\delta^*), \notag
     \end{align}
for $i=1,2$.  From straightforward modifications of Theorem 3.1 of 
\cite{AS06}, the limits
\[
   \lim_N \fN \log Z_{[0,N]}^{\{V_i\}}(\beta,u,\Phi_{iN}), \quad i=1,2,
   \]
both exist as nonrandom constants a.s.  By \eqref{E:annealmax} these
constants are at most
$(1-\theta)F(\delta^*)$, while by the free energy inequality in
\eqref{E:highfree2}, provided $\beta$ is sufficiently small (depending on $k_0,l_0$) we have 
$\beta f^q(\beta,u) >
(1-\theta)F(\delta^*)$.  This means that
\[
   \limsup_N \fN \log \mu^{\beta,u,\{V_i\}}_{[0,N]} \left( \left| \frac{L_N}{N} -
     \delta^* \right| > \lambda\delta^* \right)<0,
   \]
which establishes the contact-fraction inequality in \eqref{E:highfree2}.

We turn next to \eqref{E:lowfree1}.  Let $\eta, \nu > 0$ and let
$\{k_N\}$ be a sequence of integers with
$\chi_N = k_N/N \geq \eta$
for all $N$.  Define the events 
\[
   A_N = \left\{ \{v_i\}: \left| \sum_{i=1}^N v_i \right| \leq
     \nu^2 \beta\eta N \right\}, \quad
     A_{Nk} = A_N \cap \{L_N = k\}.
    \]
Now for large $N$, since $P^V(A_N) \geq \half$, for all
$\lambda>\beta\nu$,
\begin{align} \label{E:condlLD}
   P^V\left( \sum_{i=1}^{k_N} V_i \geq \lambda k_N\ \bigg|\ A_N \right) &\leq
     2 P^V\left( \sum_{i=1}^{k_N} V_i \geq \lambda k_N, 
     \sum_{i=k_N+1}^N V_i \leq
     -(\lambda k_N - \nu^2 \beta\eta N) \right) \notag \\
   &\leq 2 \exp\left( -k_N I_V(\lambda)
     -(N-k_N)I_V\left(-\frac{\lambda k_N - \nu^2 \beta\eta N}{N-k_N}
     \right) \right) \\
   &= 2\exp\left( -\left( I_V(\lambda) +
     \frac{1-\chi_N}{\chi_N} I_V\left( -\frac{\lambda\chi_N
     - \nu^2 \beta\eta}{1-\chi_N} \right)
     \right) k_N \right) \notag \\
   &\leq \exp\left( -\frac{\lambda^2 - 2\nu^2 \beta\eta\lambda}{2(1-\chi_N)}
     k_N \right)  \notag \\
   &\leq \exp\left( -\frac{\lambda^2(1 - 2\nu\chi_N)}{2(1-\chi_N)}
     k_N \right)  \notag \\
   &\leq \exp\left( -\frac{\lambda^2}{2(1 - (1-2\nu)\chi_N)}
     k_N \right).  \notag
\end{align}
Note that conditioning on the event $A_N$ of a ``typical disorder'' 
effectively imposes an
extra cost in the form of the second large deviation on the right 
side of the first line of
\eqref{E:condlLD}, and this extra cost is reflected in the term 
$(1-2\nu)\chi_N$ on the
right side of
\eqref{E:condlLD}.  Conditioning on $A_N$ is related to what is called 
the \emph{Morita approximation}, in which moments of the disorder are
effectively held fixed, or nearly fixed (\cite{Mo64},\cite{ORW02}.)
We have using \eqref{E:condlLD} that for large $N$,
\begin{align}
   E^V&\left( e^{\beta \sum_{i=1}^{k_N} V_i}\ \big|\ A_{Nk_N} \right) \notag \\
   &\leq
     1 + \int_0^{\beta\nu} \beta k_N e^{\beta k_N \lambda}
     P^V\left( \sum_{i=1}^{k_N} V_i \geq \lambda k_N\ \bigg|\ A_N \right)\
     d\lambda \notag \\
   &\qquad + \int_{\beta\nu}^\infty \beta k_N e^{\beta k_N \lambda}
     P^V\left( \sum_{i=1}^{k_N} V_i \geq \lambda k_N\ \bigg|\ A_N \right)\
     d\lambda \notag \\
   &\leq 1 + \beta^2 \nu k_N e^{\beta^2 \nu k_N} +
     \beta k_N \int_{-\infty}^\infty  \exp\left( \left( \beta\lambda -
     \frac{\lambda^2}{2(1 - (1-2\nu)\chi_N)} \right) k_N \right)\ 
     d\lambda \notag \\
   &\leq  \exp\left( \left( \frac{ \beta^2 (1-(1-2\nu)\chi_N) }{2}
     + \nu \right) k_N \right), \notag
\end{align}
and hence
\begin{align}
   E^V &\left( Z_{[0,N]}^{\{V_i\}}(\beta,u,\{L_N \geq \eta N\})\ 
\big|\ A_N \right)
     \notag \\
   &\leq \left\langle \exp\left( \left( \beta u + \frac{ \beta^2 (1 -
     (1-2\nu)\frac{L_N}{N}) }{2}
     + \nu \right) L_N \right) \delta_{\{L_N \geq \eta N\}} \right\rangle^X.
\end{align}
Since $\nu$ is arbitrary it follows that
\begin{align} \label{E:abovetheta}
   \limsup_N &\fN E^V\left( \log Z_{[0,N]}^{\{V_i\}}(\beta,u,\{L_N 
\geq \eta N\})\
     \big|\ A_N \right) \\
   &\leq
     \limsup_N \fN \log E^V\left( Z_{[0,N]}^{\{V_i\}}
     (\beta,u,\{L_N \geq \eta N\})\ \big|\ A_N \right) \notag \\
   &\leq \sup\left\{ \left( \beta u + \frac{ \beta^2 }{2} \right) \delta
     - \frac{\beta^2 \delta^2}{2} - \delta I_E(\delta^{-1}):
     \delta \geq \eta \right\} \notag \\
   &= \sup\left\{\beta \Delta \delta
     - \frac{\beta^2 \delta^2}{2} - \delta I_E(\delta^{-1}):
     \delta \geq \eta \right\}. \notag
\end{align}
One may view the two negative terms on the right side of 
\eqref{E:abovetheta} as two
separate costs:  the first, of order $\delta^2$, is the 
above-mentioned cost of the second ``compensating'' large deviation on the right side of the first line
of \eqref{E:condlLD}.  The second is the cost of lowering the average excursion length 
enough to get $L_N \approx \delta N$.  By \eqref{E:excurscost}, since $c>3/2$, when $\delta$ is small the cost of the compensating large deviation 
(which does not exist in the annealed model) exceeds the 
cost of lowering the average excursion length.  This is what underlies \eqref{E:lowfree1}.

For $\eta > 2\Delta/\beta$ the right side of \eqref{E:abovetheta} is at
most
\begin{equation} \label{E:negfree}
   \sup\left\{\beta \Delta \delta - \frac{\beta^2 \delta^2}{2}:
     \delta \geq \eta \right\} < 0.
     \end{equation}
As in the proof of Theorem 3.1 of \cite{AS06}, there exist constants $\beta
f^q(\beta,u,\eta^-) \geq 0$ and $\beta
f^q(\beta,u,\eta^+)$ such that
\begin{align} \label{E:thetaplus}
   \lim_N \fN \log Z_{[0,N]}^{\{V_i\}}(\beta,u,\{L_N \geq \eta N\}) &=
     \lim_N \fN E^V\left( \log Z_{[0,N]}^{\{V_i\}}(\beta,u,\{L_N \geq 
\eta N\}) \right)
     \notag \\
   &= \beta f^q(\beta,u,\eta^+) \quad \text{a.s.}
\end{align}
and
\begin{align} \label{E:thetaminus}
   \lim_N \fN \log Z_{[0,N]}^{\{V_i\}}(\beta,u,\{L_N \leq \eta N\}) &=
     \lim_N \fN E^V\left( \log Z_{[0,N]}^{\{V_i\}}(\beta,u,\{L_N \leq 
\eta N\}) \right)
     \notag \\
   &= \beta f^q(\beta,u,\eta^-) \quad \text{a.s.}
\end{align}
Further, again as in the proof of Theorem 3.1 of \cite{AS06}, by 
truncating the $V_i$ at
some large $M$ to obtain random variables $\tilde{V}_i$ and applying 
Azuma's inequality
\cite{Az67} to
\[
   \log Z_{[0,N]}^{\{\tilde{V}_i\}}(\beta,u,\{L_N \geq \eta N\}),
\]
since $\liminf_N P^V(A_N)>0$ we obtain using \eqref{E:thetaplus} and \eqref{E:thetaminus} that
\begin{equation} \label{E:azuma}
   \left| E^V\left( \log Z_{[0,N]}^{\{V_i\}}(\beta,u,\{L_N \geq \eta N\})\ \big|\ A_N
     \right) - E^V\left( \log Z_{[0,N]}^{\{V_i\}}(\beta,u,\{L_N \geq \eta N\}) 
     \right) \right| = o(N).
\end{equation}
With \eqref{E:abovetheta} and \eqref{E:thetaminus} this shows that $\beta
f^q(\beta,u,\eta^+)<0$.  We can then conclude from \eqref{E:thetaplus} and
\eqref{E:thetaminus} that
\[
   \mu^{\beta,u,\{V_i\}}_{[0,N]}(L_N \geq \eta N) \to 0 \quad \text{ 
as } N \to \infty,
\]
proving the contact-fraction inequality in \eqref{E:lowfree1}.

It follows from \eqref{E:abovetheta} and \eqref{E:azuma} (with $\eta=0$) that
\[
   \beta f^q\left( \beta, -\frac{\beta}{2} + \Delta \right) \leq
     \sup\left\{\beta \Delta \delta - \frac{\beta^2 \delta^2}{2}:
     \delta \geq 0 \right\} = \frac{\Delta^2}{2},
\]
which is the free-energy inequality in \eqref{E:lowfree1}.

\section{Proof of Theorems \ref{T:smallc}, \ref{T:c32} and \ref{T:transient}} \label{S:moreproof}

\begin{proof}[Proof of Theorem \ref{T:c32}] For the free-energy inequality, the only changes needed from the case $3/2 < c < 2$ in Theorem \ref{T:bigc} involve the definition of $p_M$ and the fact that for $c=3/2$, \eqref{E:Deltalower} is not a sufficient condition for the first inequality in \eqref{E:pMlower}.  From Lemma \ref{L:boundreturns} the proper choice is now $p_M = 1/K_{53}\tphi(M)$, and for some $K_{57} = K_{57}(\epsilon)$  the first inequality in \eqref{E:pMlower} then holds provided 
\[
  \tphi(M) \leq \frac{K_{57}}{\beta^2},
  \]
for which, by \eqref{E:alpha0val}, it suffices that
\[
  \frac{K_{58} \ophi_{1/2}^*\left( \frac{1}{\beta\Delta} \right)^2 }{ (\beta\Delta)^2 } =
    \frac{K_{58}}{ (\beta\Delta)^2 \hat{\varphi}_{1/2}\left( \frac{1}{\beta\Delta} \right) } 
     \leq \tphi^{-1}\left( \frac{K_{57}}{\beta^2} \right),
    \]
or 
\[
  \frac{K_{59}}{\beta\Delta} \leq \tphi^{-1}\left( \frac{K_{57}}{\beta^2} \right)^{1/2}
    \ophi_{1/2}\left( \tphi^{-1}\left( \frac{K_{57}}{\beta^2} \right)^{1/2} \right)
    = \frac{ \tphi^{-1}\left( \frac{K_{57}}{\beta^2} \right)^{1/2} }
    { \varphi\left( \tphi^{-1}\left( \frac{K_{57}}{\beta^2} \right) \right) },
    \]
which is equivalent to $\Delta \geq \Delta_0$, for an appropriate choice of $K_{5}$.  Here all $K_i$ depend on $\epsilon$.

The proof of the contact-fraction inequality in \eqref{E:highfree1} from the free-energy inequality remains unchanged from Theorem \ref{T:bigc}.
\end{proof}

\begin{proof}[Proof of Theorem \ref{T:smallc}]
Most of the proof of the free-energy inequality in \eqref{E:highfree1} is the same as that of the free-energy inequality in \eqref{E:highfree2}, but we have the following changes.
As in the proof of \eqref{E:highfree2}, let $\{Y_i\}$ be an independent copy of the Markov chain $\{X_i\}$, under the distribution $P^X$.  Under the hypotheses of the theorem, it follows from Lemma \ref{L:Erickson} that $\{(X_i,Y_i)\}$ is transient.  This means that $B_M(\{X_n\},\{Y_n\}) + 1$ is stochastically smaller than a geometric random variable with parameter $\tilde{p} = P^X(\tE_1 = \infty)$, where we recall that $\tE_1$ is the length of the first excursion for the chain $\{(X_i,Y_i)\}$.  Thus the dependence of $p_M$ on $M$ is effectively removed--we can achieve \eqref{E:pMlower} (with $p_M$ replaced by $\tilde{p}$) merely by taking $\beta$ sufficiently small, and \eqref{E:Deltalower} is not needed.

The proof of the contact-fraction inequality in \eqref{E:highfree1} from the free-energy inequality remains unchanged from Theorem \ref{T:bigc}.
\end{proof}

\begin{proof}[Proof of Theorem \ref{T:transient}]
We have
\[
  Z^{\{V_i\}}_{[0,N]}(\beta,u) = \left\langle e^{\beta u_c^d L_N} \right\rangle^X
    \sum_{\{x_i\}} \exp\left(\beta\sum_{i=1}^N (u - u_c^d + V_i) 
    \delta_{\{x_i = 0\}} \right) \mu_{[0,N]}^{\beta,u_c^d,0}(x_{[0,N]})
    \]
and by definition of $u_c^d$ and continuity of the free energy,
\[
  \lim_N \fN \log \left\langle e^{\beta u_c^d L_N} \right\rangle^X = 0.
  \]
Hence to show the free energies are equal, it suffices to show that there exist constants $a_N$ with $\log a_N = o(N)$ and 
\begin{equation} \label{E:sufficient}
  \frac{1}{a_N} \leq \frac{ d\mu_{[0,N]}^{\beta,u_c^d,0} }{ dP^{X,R}_{[0,N]} }(\{x_i\}) \leq a_N
    \quad \text{for all } \{x_i\}.
    \end{equation}
Let $T_N$ denote the time of the last return to 0 in $[0,N]$, and let $\oF(x) = P^X(x < E_1 < \infty)$,
so $\oF(x) \sim (c-1)^{-1}x^{-(c-1)}\varphi(x)$.
The main observation is that by \eqref{E:ucd}, the difference between $\mu_{[0,N]}^{\beta,u_c^d,0}$ and $P^{X,R}_{[0,N]}$ involves only the final excursion in progress at time $N$, in the sense that
\[
  \frac{ d\mu_{[0,N]}^{\beta,u_c^d,0} }{ dP^{X,R}_{[0,N]} } = 
    q_N \frac{ \oF(N-T_N) + P^X(E_1 = \infty) }{ \oF(N-T_N) }    
    \]
where 
\[
  \frac{1}{q_N} = \left\langle \frac{ \oF(N-T_N) + P^X(E_1 = \infty) }{ \oF(N-T_N) } \right\rangle^{X,R}.
  \]
Here $\left\langle \cdot \right\rangle^{X,R}$ denotes expectation with respect to $P^{X,R}$.  We have
\[
  \frac{1}{q_N} \geq K_{60}\left\langle (N-T_N)^{c-1} \varphi(N-T_N)^{-1} \right\rangle^{X,R}
    \geq K_{61} N^{c-1} \varphi(N)^{-1} P^{X,R}\left( T_N \leq \frac{N}{2} \right),
  \]
and as in \eqref{E:CNcondl} and the calculations following it we have (since $c<2$)
\[
  P^{X,R}\left( T_N \leq \frac{N}{2} \right) \geq K_{62} \quad \text{for all } N.
  \]
Therefore
\[
  \frac{ d\mu_{[0,N]}^{\beta,u_c^d,0} }{ dP^{X,R}_{[0,N]} } \leq \frac{q_N}{\oF(N)} \leq K_{63}.
  \]
In the other direction,
\[
  \frac{1}{q_N} \leq \left\langle \frac{1}{\oF(N - T_N)} \right\rangle^{X,R}
    \leq K_{64} \left\langle \frac{(N-T_N)^{c-1}}{ \varphi(N-T_N) } \right\rangle^{X,R} 
    \leq K_{65} \frac{N^{c-1}}{\varphi(N)},
  \]
so
\[
  \frac{ d\mu_{[0,N]}^{\beta,u_c^d,0} }{ dP^{X,R}_{[0,N]} } \geq 
    K_{66} \frac{\varphi(N)}{N^{c-1}},
  \]
so \eqref{E:sufficient} is proved.

Equality of the contact fractions follows immediately from equality of the free energies, by definition of the contact fraction.
\end{proof}

\begin{section}{The Excluded Cases $c=1$ and $c \geq 2$} \label{S:excludedc}
We conclude with a few remarks about the cases $c=1$ and $c \geq 2$ not covered by our results.  

For $c=1$, the left side of \eqref{E:fracbound} is not bounded below uniformly in $\Delta$, so a new proof of (P1) is needed.  More importantly, in this case there are two disagreeing definitions of the correlation length.  One, which we denote here by $M^*$, is given by \eqref{E:Mdef} (i.e. $M^*$ is the $M$ we used above), and the other, which we denote by $M_0$, is suggested by \eqref{E:ma0}:  $M_0 = K_{9}/\alpha_0$, the inverse of the free energy.  For $c>1$ these two values are (by \eqref{E:ma0}) asymptotically the same as $\Delta \to 0$, but for $c=1$ we have $M_0/M^* \to \infty$ as $\beta\Delta \to 0$.  Inequality (P2) can be interpreted as saying that the entropy cost of the event $\Xi_N$ is at most a small mulltiple of $1/M^*$ per unit length, but for $c=1$ a small multiple of $1/M^*$ becomes an unacceptably large multiple of $1/M_0$.  Thus for $c=1$ something else must substitute for the event $\Xi_N$, which plays the role of creating an approximate renewal at the start of each block of length $N$, leading to the factoring of the partition function expressed in (P1).

The case $c \geq 2$ is really two cases:  $c=2$ with $\langle E_1 \rangle^X = \infty$, which means the transition in the annealed system is continuous, and $\langle E_1 \rangle^X < \infty$ which means the transition is discontinuous.  For $c<2$, the idea of Lemma \ref{L:doubletail} is that $\tL_n = k$ typically means the $k$th excursion was one of the first few excursions of length at least of order $n$, so $\tL_n$ is of the same order as a geometric random variable with parameter $1/P^X(\tE_1 > n)$.  This will not be valid for $c \geq 2$, so a different approach is required.  Overall, though, despite significant changes in the details, the core ideas of our proof should carry over well to the case of $c=2$ with $\langle E_1 \rangle^X = \infty$.  Recall that $\Delta_0(\beta)$ is (up to a constant) the magnitude of $\Delta$ for which the upper bound $2\Delta/\beta$ is equal to the annealed contact fraction.  When the transition is discontinuous (i.e. when $\langle E_1 \rangle^X < \infty$), we have $C^a(\beta,u) > 1/\langle E_1 \rangle^X$ for all $u>u_c^a(\beta)$, and therefore $\Delta_0(\beta)$ is not $o(\beta)$ (and hence not $o(u_c^a(\beta)$) as $\beta \to 0$.  Thus for $c>2$ it is not accurate to describe the disorder-induced changes as being confined to a very small interval above $u_c^a$.
\end{section}

\end{document}